\documentclass[a4paper, 11pt]{article}

\usepackage[margin=1.2in]{geometry}
\usepackage{parskip}
\usepackage[titletoc,title]{appendix}

\usepackage{array}
\usepackage{caption}
\usepackage{rotating}
\usepackage[hidelinks, colorlinks = true, linkcolor = blue, urlcolor = blue, anchorcolor = blue, citecolor = blue]{hyperref}
\usepackage{blkarray, bigstrut}
\usepackage{multirow}
\usepackage{booktabs}
\usepackage{siunitx}
\usepackage{mathtools}
\usepackage{amsmath}
\usepackage{amssymb}
\usepackage{amsthm}
\usepackage{complexity}
\usepackage{natbib}
\usepackage{graphicx}
\usepackage{pgf}
\usepackage{tikz}
\usetikzlibrary{matrix,patterns,patterns.meta,positioning,decorations,arrows,shapes,decorations.markings,calc}

\usetikzlibrary{arrows.meta,bending}
\usepackage[latin1]{inputenc}
\usepackage{verbatim}
\usepackage{subfiles}
\usepackage{pgfplots}
\pgfplotsset{compat=newest}
\usepgfplotslibrary{groupplots}
\usepackage{color}
\usepackage{float}
\usepackage{enumerate}
\usepackage{emptypage}
\usepackage{subcaption}
\usepackage{setspace}
\usepackage[export]{adjustbox}
\usepackage{subfiles}
\usepackage[algoruled]{algorithm2e}
\usepackage{xcolor}
\newcolumntype{H}{>{\setbox0=\hbox\bgroup}c<{\egroup}@{}}

\usepackage{pdflscape}

\newcommand{\mc}{\mathcal}

\usepackage[normalem]{ulem}

\usepackage[short]{optidef}
\usepackage{todonotes}
\setuptodonotes{inline}

\title{Line Planning at Scale: Models, Methods, and Insights}

\author{R.N. van Lieshout$^*$, B.T.C. van Rossum,
\vspace{0.1cm}\\
\small{Department of Industrial Engineering \& Innovation Sciences} \\
\small{Eindhoven University of Technology, The Netherlands}
\vspace{0.1cm}\\
\small{$^\star$Corresponding author. Email: r.n.v.lieshout@tue.nl}
}

\begin{document}

%title
\maketitle

%abstract
\begin{abstract}
\noindent 
Line planning, the problem of deciding which lines to operate and at what frequency, is a fundamental step in public transport planning. To accurately model passenger routing, the problem is traditionally defined on a change-and-go network (CGN), which captures transfers between lines exactly. However, this network grows large quickly and is hard to solve at scale. We compare the CGN against three more compact models, differing with respect to how transfers are approximated, and characterize how they relate in terms of solution quality and modeling accuracy. We develop state-of-the-art solution methods tailored to each model, and evaluate all four across 972 instances based on the Dutch and Swiss railway networks. Contrary to the CGN's canonical status, we find that it is competitive only on small or easy instances and often fails to find any feasible solution on large networks. Instead, a compact direct connection model performs best overall, finding the best solution on over 83\% of instances. Our results indicate that carefully designed approximations, rather than exact transfer modeling, are the more promising foundation for large-scale line planning.
\vspace{5mm}
\end{abstract}

\section{Introduction}
\label{sec:intro}

Line planning is a fundamental step in public transit planning, concerned with determining which lines to operate and at what frequency. A good line plan must balance the needs of passengers --- minimizing travel times and transfers --- against the operational costs of the operator, subject to infrastructure and capacity constraints. Despite its conceptual simplicity, line planning gives rise to challenging mixed-integer programs, especially when passenger routing and capacity constraints are modeled explicitly.

A central difficulty is how to represent passenger routes and transfers within the optimization model. The \textit{change-and-go network} (CGN) of \citet{schobel2006line} provides the canonical solution: by tracking which line each passenger is traveling on at every point in their journey, it enables exact modeling of transfers and per-line capacity constraints. As such, it has served as the foundation for most subsequent work in line planning with passenger routing. However, the CGN grows rapidly as the size of the line pool increases. Despite the use of column generation, solving even the linear programming (LP) relaxation of large instances remains challenging. 

These difficulties have motivated a range of alternative representations that approximate transfers with varying degrees of accuracy. At one extreme, \citet{borndorfer2007column} route passengers directly over the public transport network (PTN), ignoring transfers entirely. While this yields a far more compact formulation, neglecting transfer information can lead to poor solutions when transfer penalties are significant. Two more recent approaches seek a middle ground, recovering part of this transfer information without incurring the CGN's computational cost. \citet{borndorfer2015metric} introduce metric inequalities within the PTN framework to distinguish passengers on a direct route from those requiring a transfer. \citet{yao2024railway} instead propose a \textit{direction-expanded PTN} (DPTN) that tracks the direction of travel rather than the line, capturing transfers more accurately than the PTN while remaining more compact than the CGN. Despite their appeal, the theoretical relationships among these approaches and their relative computational performance have not been systematically studied.

This paper provides a unified theoretical and computational treatment of the various modeling approaches for passenger-oriented line planning. Our central finding overturns a common assumption: although the change-and-go network models transfers exactly, it is rarely the best choice for large-scale line planning. Because it grows so large, it frequently fails to produce even a feasible solution on realistic instances, whereas carefully designed approximations that partially capture transfers consistently yield better line plans. The CGN is thus perhaps best understood not as the default optimization model but as a benchmark that is invaluable for evaluating line plans, yet seldom the best tool for producing them. This motivates a closer study of the approximate models that outperform it.

We arrive at this finding by first recasting the CGN, the PTN, the direct connection model of \citet{borndorfer2015metric}, and the direction-expanded network of \citet{yao2024railway} within a common modeling framework, making their structural relationships explicit. We establish a model hierarchy in terms of solution quality and show that the direct connection and direction-expanded models are not directly comparable, each dominating the other
on some instances. We further show that the exact direct connection constraints
of \citet{borndorfer2015metric} can be reinterpreted as Benders cuts, which provides a transparent derivation of the metric inequalities and suggests natural extensions.

To compare the models in practice, we develop state-of-the-art solution approaches tailored to each network structure, including a branch-price-and-cut framework and a Benders decomposition. We then conduct an extensive computational study on 972 instances based on real railway networks in the Netherlands and Switzerland, varying the problem parameters of practical interest such as line pool size, budget, and transfer penalty. The study yields concrete and sometimes counterintuitive guidance. The direct connection model performs best overall, finding the best solution on over 83\% of instances, while the exact CGN is competitive only on small or easy instances and fails to find any feasible solution on more than half of the full set. The direction-expanded model scales well, with its relative performance strongest on the largest instances and on large line pools. These findings suggest that the additional complexity of exact transfer modeling is rarely warranted. Instead, approximate but compact
representations are the more promising foundation for large-scale line planning. To support future work, we make all code, instances, and best-known solutions publicly available at \url{https://github.com/rn-van-lieshout/line-planning-benchmark}.
 %[link removed to preserve anonymity] 
 
The remainder of the paper is organized as follows. Section~\ref{sec:literature} reviews the relevant literature. Section~\ref{sec:problem} formulates the line planning problem and introduces the CGN. Section~\ref{sec:models} presents the PTN-based and DPTN models and establishes their theoretical relationships. Section~\ref{sec:sol} describes the solution approaches and acceleration strategies. Section~\ref{sec:experiments} reports the results of our computational study, and Section~\ref{sec:con} concludes.
\section{Literature Review}
\label{sec:literature}

Line planning is one of several interdependent planning stages in public transport planning, which also includes network design, timetabling, and rolling
stock and crew scheduling. We refer to \citet{schmidt2024planning} for a recent overview of the line planning literature, and restrict our discussion here to contributions considering line planning with passenger routing.

Two early papers have been particularly influential in shaping the line planning literature. \citet{schobel2006line} introduce the change-and-go network as a modeling framework for passenger routing with explicit
transfer penalties, providing the formal foundation on which much subsequent work is built. The exactness of this representation comes at a computational cost: despite the use of column generation, the authors struggle to solve the LP relaxation of large instances, observing that ``the main problem of our approach is the size of the change\&go network depending mainly on the size of the line pool". \citet{borndorfer2007column} prevent this curse of dimensionality by routing passengers on the public transport network, effectively ignoring transfers. Although the compactness of the PTN allows passenger routes to be
enumerated in advance, this paper is widely credited with establishing column generation as the method of choice for large-scale line planning, with routes generated dynamically.

Several subsequent papers build directly on the CGN, with extensions addressing user-optimal route choice \citep{GOERIGK2017424}, partial integration with timetabling \citep{BURGGRAEVE2017134}, frequency-dependent transfer times \citep{bull2019optimising}, mode choice \citep{hartleb2023modeling}, multi-period planning \Citep{van2025multi}, and crowding \citep{lu2025line}. Despite considerable computational advances, these papers either remain limited to medium-scale instances or resort to heuristics -- such as genetic algorithms or column generation-based heuristics -- to handle larger ones. In a cost-oriented setting, \citet{friedrich2017integrating} use the CGN not as an optimization model but as a routing network for a heuristic passenger assignment, projecting the resulting flows onto the PTN to feed a cost-minimizing line model. They find CGN-based assignment yields better line concepts than PTN-based assignment, but do not compare the two as exact formulations.

Only three papers attempt to circumvent the CGN's size bottleneck through the use of the PTN or other smaller network representations. \citet{borndorfer2015metric}
propose routing passengers over the PTN rather than the CGN, and introduce metric inequalities to recover the ability to charge transfer penalties --- an approach we refer to as the PTN-based direct connection model (PTN-DC). Despite promising
computational results, this approach went largely unnoticed by the line planning community until \citet{GATT2025100164} recently extended it with service-level
constraints, solving the resulting model via a heuristic combining column generation and column enumeration. \citet{yao2024railway} take a different approach, proposing a direction-expanded representation --- which we refer to as the DPTN ---  that is more compact than the CGN while modeling transfers more accurately than the PTN. Since transfers are still only partially captured,
they embed the model in a two-phase approach in which the second phase solves a CGN-like model with exact transfers over the neighborhood of the first-phase line plan.

\begin{table}[h]
\centering
\caption{Overview of modeling approaches, instance sizes, and solution methods in the line planning literature.}
\label{tab:litreview}
\small
\resizebox{\textwidth}{!}{
\setlength{\tabcolsep}{4pt}
\begin{tabular}{@{}llrrl p{4.3cm}@{}}
\toprule
Paper & Model & Pool & ODs & Method & Notes \\
\midrule
\citet{schobel2006line}   & CGN & --   & --       & CG            & LP relaxation only,\ no integer solution \\
\citet{GOERIGK2017424}    & CGN & 132  & 48{,}842 & GA            & User-optimal route choice; MIP for small instances \\
\citet{BURGGRAEVE2017134} & CGN & --   & 4{,}645  & MIP           & Freq.-dep.\ transfer times; partial TT integration \\
\citet{bull2019optimising}& CGN & 174  & 4{,}645  & LP heur.\     & Freq.-dep.\ transfer penalties \\
\citet{hartleb2023modeling}& CGN& 64   & 174      & MIP           & Mode choice \\
\Citet{van2025multi}      & CGN & 10   & 1{,}406  & MIP           & Multi-period; variable stopping patterns \\
\citet{lu2025line}        & CGN & 534  & 56{,}916 & CG heur.\     & Crowding \\
\addlinespace
\citet{borndorfer2007column} & PTN & 80 & 4{,}685  & CG heur.\     & Lines priced in some experiments \\
\citet{friedrich2017integrating} & PTN & 275 & 567 & MIP        & CGN used as heuristic routing network \\
\addlinespace
\citet{borndorfer2015metric} & PTN-DC & 4{,}342 & 7{,}734 & Metric ineq.\ & \\
\citet{GATT2025100164}    & PTN-DC & LG & 528      & CG heur.\     & Service-level constraints \\
\addlinespace
\citet{yao2024railway}    & DPTN & LG  & 934      & B\&P + LS\     & No fixed line cost \\
\bottomrule
\end{tabular}
}
\smallskip
{\footnotesize
\textit{Abbreviations:} LG = line generation (no fixed pool); CG = column generation;
GA = genetic algorithm; B\&P = branch-and-price; LS = local search; heur.\ = heuristic; TT = timetabling.}
\end{table}

Table~\ref{tab:litreview} summarizes the modeling approaches, instance size, and solution method of the works most closely related to ours. A few patterns stand out. The CGN underpins the majority of the literature, typically solved exactly via mixed-integer programming (MIP) or column generation but only on small to medium instances, while the largest instances are handled either by the compact PTN or by heuristics. The PTN-DC and DPTN representations, despite their promise, have received comparatively little attention: only two papers study the former and just one the latter. 

The present paper provides the first systematic comparison of these alternative models. We formalize the PTN-based and direction-expanded approaches within a common modeling framework, establishing theoretical relationships between the resulting formulations. We develop tailored state-of-the-art solution methods to solve the formulations, and evaluate their computational trade-offs across a benchmark of 972 instances. These instances lie at the upper end of the spectrum in Table~\ref{tab:litreview}: with line pools of up to 12{,}359 candidate lines and up to 25{,}510 OD pairs, we exceed all
previously reported line pools and rank among the largest in OD count, with the few studies reporting more relying on heuristics.
\section{Problem Formulation}
\label{sec:problem}

The line planning problem asks to select a set of lines and assign passengers to routes so as to minimize total passenger travel time, subject to capacity, budget, and frequency constraints. The starting point for a formal definition is the \textit{public transport network}, an undirected graph $G^\text{PTN} = (\mc{S}, \mathcal{C})$ containing a node
for each station and an edge for each pair of stations that can be connected by a direct (non-stop) service. Each \textit{connection} $c \in \mathcal{C}$ has a travel time $\tau_c$ and a maximum service frequency $f_c$. For metro or tram networks, the PTN reflects the physical infrastructure, with connections corresponding to tracks. In urban or national rail networks, express services often skip intermediate stations, so that non-adjacent stations may be directly connected in the PTN.

A \textit{line} $\ell$ defines a simple directed path in $G^\text{PTN}$. The set of candidate lines, called the \textit{line pool}, is denoted by $\mathcal{L}$. Infrastructural constraints --- for example, whether track switches allow a train to enter a station from one direction and exit in another --- can be accounted for during the construction of $\mathcal{L}$. Each line $\ell$ can be operated at a non-negative integer frequency, incurring a variable cost $c^V_\ell$ per unit frequency and a fixed, frequency-independent cost $c^F_\ell$. The line provides a passenger capacity of $\kappa_\ell$ per unit frequency. For clarity of exposition, we treat lines as directed. The extension to bidirectional lines, where each line is operated in both directions, is straightforward. The total budget available is denoted by $B$.

Passengers are represented as a set of origin-destination (OD) pairs $\mathcal{OD} \subseteq \mathcal{S} \times \mathcal{S}$, indexed by $k$, where $v_k$ denotes the demand volume of pair $k$. To ensure feasibility, we allow the operator to leave a portion of demand unserved, at a penalty cost of $\pi_k$ per unserved passenger for OD pair $k$. Served passengers are assigned to a route from their origin to their destination, whose travel time is the sum of the travel times of all traversed connections, the dwell times at intermediate stations, and a penalty $\sigma$ for each transfer between lines. Capacity constraints must be satisfied on every arc used by passengers. 

These considerations motivate the use of a line-expanded network --- the \textit{change-and-go network} --- tracking which lines passengers board and where they transfer. Formally, the CGN is a directed graph $G^\text{CGN} = (\mathcal{V}, \mathcal{A})$. The node set $\mathcal{V}$ contains one \textit{station node} for each station $s \in \mathcal{S}$, used for boarding, alighting or transferring, and one \textit{travel node} $(\ell, s)$ for each line $\ell \in \mathcal{L}$ and station $s \in \mathcal{S}$ visited by $\ell$, used for traveling between stations. The arc set $\mathcal{A}$ consists of two types:
\begin{itemize}
    \item \textit{Running arcs} $\mathcal{A}^R$: an arc from $(\ell, s)$ to
        $(\ell, s')$ for each consecutive station pair $(s, s')$ on line
        $\ell$, representing travel along a connection.
    \item \textit{Transfer arcs} $\mathcal{A}^T$: arcs between the station
        node of $s$ and each travel node $(\ell, s)$ at $s$, representing
        boarding, alighting, or transferring between lines.
\end{itemize}
The travel time of arc $a \in \mathcal{A}$, capturing either running time or a transfer penalty, is denoted $\tau_a$. Dwell times may be incorporated into running times, with a corresponding correction applied to transfer arc durations. For a running arc $a$, we write $\ell(a) \in \mathcal{L}$ for its
associated line. Figure~\ref{fig:ptn-cgn} illustrates a small PTN with two lines and the corresponding CGN. 

% Required packages:
% \usepackage{tikz}
% \usetikzlibrary{arrows.meta, bending}
% \usepackage{xcolor}
% \usepackage{subcaption}

\definecolor{tol_blue}{HTML}{4477AA}
\definecolor{tol_red}{HTML}{EE6677}

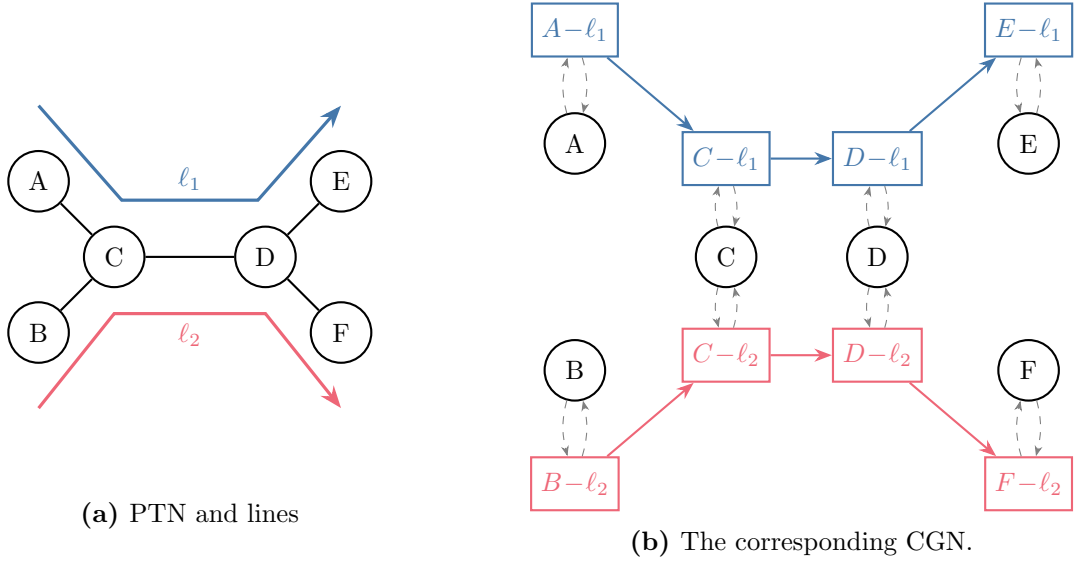
\begin{figure}[htbp]
\centering

\begin{subfigure}[c]{0.35\linewidth}
\centering
% =====================================================================
%  (a) PTN
% =====================================================================
\begin{tikzpicture}[
  stn/.style   = {circle, draw, thick, minimum size=8mm, font=\small},
  conn/.style  = {thick},
  line1/.style = {very thick, tol_blue},
  line2/.style = {very thick, tol_red},
]
  % Phantom extents: match CGN vertical span (-3..3) and center on y=0,
  % so that the y=0 line (C and D) aligns across both subfigures.
  \path (0,-3.0) (0,3.0);

  \node[stn] (A) at (-2, 1) {A};
  \node[stn] (B) at (-2,-1) {B};
  \node[stn] (C) at (-1, 0) {C};
  \node[stn] (D) at ( 1, 0) {D};
  \node[stn] (E) at ( 2, 1) {E};
  \node[stn] (F) at ( 2,-1) {F};

  \draw[conn] (A) -- (C);
  \draw[conn] (B) -- (C);
  \draw[conn] (C) -- (D);
  \draw[conn] (D) -- (E);
  \draw[conn] (D) -- (F);

  % Line 1: A->C->D->E
  \draw[line1, -{Stealth}]
    (-2, 2.0) -- (-0.9, 0.75) -- (0.9, 0.75) -- (2, 2.0);
  \node[tol_blue, font=\small, above] at (0, 0.75) {$\ell_1$};
  \draw[line2, -{Stealth}]
    (-2,-2) -- (-1,-0.75) -- (1,-0.75) -- (2,-2);
  \node[tol_red, font=\small, below] at (0,-0.75) {$\ell_2$};

\end{tikzpicture}
\caption{PTN and lines}
\label{fig:ptn}
\end{subfigure}
\hspace{1.5cm}
\begin{subfigure}[c]{0.50\linewidth}
\centering
% =====================================================================
%  (b) CGN
% =====================================================================
\begin{tikzpicture}[
  stn/.style   = {circle, draw, thick, minimum size=8mm, font=\small},
  trv1/.style  = {rectangle, draw, thick, tol_blue, minimum width=1.1cm,
                  minimum height=7mm, font=\small},
  trv2/.style  = {rectangle, draw, thick, tol_red, minimum width=1.1cm,
                  minimum height=7mm, font=\small},
  run1/.style  = {-{Stealth}, thick, tol_blue},
  run2/.style  = {-{Stealth}, thick, tol_red},
  trf/.style   = {-{Stealth[length=4pt]}, gray, thin, dashed},
]
  % Phantom extents: match PTN, center on y=0 (CGN already spans -3..3).
  \path (0,-3.0) (0,3.0);

  % Travel nodes line 1 (blue)
  \node[trv1] (L1A) at (-3,  3.0) {$A\!-\!\ell_1$};
  \node[trv1] (L1C) at (-1,  1.3) {$C\!-\!\ell_1$};
  \node[trv1] (L1D) at ( 1,  1.3) {$D\!-\!\ell_1$};
  \node[trv1] (L1E) at ( 3,  3.0) {$E\!-\!\ell_1$};

  % Travel nodes line 2 (red)
  \node[trv2] (L2B) at (-3, -3.0) {$B\!-\!\ell_2$};
  \node[trv2] (L2C) at (-1, -1.3) {$C\!-\!\ell_2$};
  \node[trv2] (L2D) at ( 1, -1.3) {$D\!-\!\ell_2$};
  \node[trv2] (L2F) at ( 3, -3.0) {$F\!-\!\ell_2$};

  % Station nodes
  \node[stn] (sA) at (-3,  1.5) {A};
  \node[stn] (sB) at (-3, -1.5) {B};
  \node[stn] (sC) at (-1,  0)   {C};
  \node[stn] (sD) at ( 1,  0)   {D};
  \node[stn] (sE) at ( 3,  1.5) {E};
  \node[stn] (sF) at ( 3, -1.5) {F};

  % Running arcs
  \draw[run1] (L1A) -- (L1C);
  \draw[run1] (L1C) -- (L1D);
  \draw[run1] (L1D) -- (L1E);

  \draw[run2] (L2B) -- (L2C);
  \draw[run2] (L2C) -- (L2D);
  \draw[run2] (L2D) -- (L2F);

  % Transfer arcs
  \draw[trf, bend left=15]  (sA)  to (L1A);
  \draw[trf, bend left=15]  (L1A) to (sA);

  \draw[trf, bend right=15] (sB)  to (L2B);
  \draw[trf, bend right=15] (L2B) to (sB);

  \draw[trf, bend left=15]  (sC)  to (L1C);
  \draw[trf, bend left=15]  (L1C) to (sC);
  \draw[trf, bend right=15] (sC)  to (L2C);
  \draw[trf, bend right=15] (L2C) to (sC);

  \draw[trf, bend left=15]  (sD)  to (L1D);
  \draw[trf, bend left=15]  (L1D) to (sD);
  \draw[trf, bend right=15] (sD)  to (L2D);
  \draw[trf, bend right=15] (L2D) to (sD);

  \draw[trf, bend right=15] (sE)  to (L1E);
  \draw[trf, bend right=15] (L1E) to (sE);

  \draw[trf, bend left=15]  (sF)  to (L2F);
  \draw[trf, bend left=15]  (L2F) to (sF);

\end{tikzpicture}
\caption{The corresponding CGN.}
\label{fig:cgn}
\end{subfigure}

\caption{The public transportation network (PTN) with two lines $\ell_1$ and $\ell_2$, and the corresponding change-and-go network (CGN). In the CGN, solid arcs are
running arcs and dashed arcs are transfer arcs; rectangular nodes represent
travel nodes $(\ell, s)$ and circular nodes represent station nodes $s$.}
\label{fig:ptn-cgn}
\end{figure}

The CGN induces a set of passenger routes $\mathcal{P}$, defined as paths in $G^{\text{CGN}}$ connecting station nodes. We let $\mathcal{P}_k \subseteq \mathcal{P}$ denote the set of routes for OD pair $k$, $\mathcal{P}(a) \subseteq \mathcal{P}$ the set of routes using arc $a$, and $\tau_p$ the
total travel time of route $p$.

To define the mathematical model, we introduce four types of decision variables. Integer variables $y_\ell \in \mathbb{Z}_{\ge 0}$ denote the operating frequency of line $\ell \in \mathcal{L}$, and binary variables $z_\ell \in \{0, 1\}$ indicate whether line $\ell$ is operated at all. The latter are required to correctly model fixed line costs. Continuous variables $x_p \ge 0$ denote the passenger flow on route $p \in \mathcal{P}$, and continuous variables $u_k \ge 0$ denote the unserved demand for OD pair $k \in \mathcal{OD}$. We are now ready to state the line planning problem, which we label \texttt{CGN} after the network representation it is built on:
\begin{mini!}
{}
{ \sum_{p \in \mathcal{P}} \tau_p x_p +
    \sum_{k \in \mathcal{OD}} \pi_k u_k }
{}{\texttt{CGN}=}
\addConstraint{ \sum_{p \in \mathcal{P}_k} x_p + u_k }{ \ge v_k }{
    \forall k \in \mathcal{OD} \label{eq:demand}}
\addConstraint{ \sum_{p \in \mathcal{P}(a)} x_p }{ \le \kappa_{\ell(a)} y_{\ell(a)}
    \quad }{ \forall a\in \mc{A}^R
    \label{eq:capacity}}
\addConstraint{ \sum_{\ell \in \mathcal{L}(c)} y_\ell }{ \le f_c \quad }{
    \forall c \in \mathcal{C} \label{eq:frequency}}
\addConstraint{ y_\ell }{ \le \bar{f}_\ell z_\ell \quad }{
    \forall \ell \in \mathcal{L} \label{eq:linking}}
\addConstraint{ \sum_{\ell \in \mathcal{L}} \left( c^V_\ell y_\ell +
    c^F_\ell z_\ell \right) }{ \le B \label{eq:budget}}
\addConstraint{ u_k }{ \ge 0 }{ \forall k \in \mathcal{OD} }
\addConstraint{ x_p }{ \ge 0 }{ \forall p \in \mathcal{P} }
\addConstraint{ y_\ell }{ \in \mathbb{Z}_{\ge 0} }{ \forall \ell \in
    \mathcal{L} }
\addConstraint{ z_\ell }{ \in \{0, 1\} }{ \forall \ell \in
    \mathcal{L}. }
\end{mini!}
The objective minimizes total passenger travel time and the
penalty for unserved demand. Constraints~\eqref{eq:demand} ensure that each OD pair is either served or accounted for as unserved demand. Constraints~\eqref{eq:capacity} enforce that the passenger flow on each arc of line $\ell$ does not exceed its capacity $\kappa_\ell y_\ell$. Constraints~\eqref{eq:frequency} limit the total frequency on each connection. Constraints~\eqref{eq:linking} link the frequency and activation variables, ensuring $z_\ell = 1$ whenever line $\ell$ is operated, where $\bar{f}_\ell = \min_{c \in \ell} f_c$ is a natural upper bound on the frequency of line $\ell$. Finally, Constraint~\eqref{eq:budget} imposes a
budget on total operating costs.

The model can alternatively be formulated with binary line-frequency variables, which avoids undesirable frequencies and enables frequency-dependent transfer times, at the cost of a significantly larger model. We do not pursue this variant here.
\section{Models}
\label{sec:models}

Section~\ref{sec:problem} presented the line planning problem using the change-and-go network, which models passenger routing and transfers exactly. However, directly solving this model is often computationally intractable for realistic instances. The CGN contains a node for every (station, line) pair, and a path for every possible sequence of lines a
passenger might use. For networks with hundreds of stations and a large line pool, the number of passenger routes $|\mathcal{P}|$ becomes enormous, making even the LP relaxation difficult to solve directly. This motivates the use of a branch-and-price approach, in which routes are generated dynamically by solving a pricing problem over the CGN. Even so, the size of the CGN itself can be a bottleneck, as the pricing problem becomes harder as the network grows.

These computational challenges have motivated the development of alternative, more compact network representations that trade modeling accuracy for tractability. In this section, we discuss two such alternatives: models based on the PTN (Section~\ref{sec:sn}) and the direction-expanded PTN
(Section~\ref{sec:desn}). All models share the same line frequency variables
$y_\ell \in \mathbb{Z}_{\ge 0}$ and binary activation variables $z_\ell \in
\{0,1\}$ introduced in Section~\ref{sec:problem}, along with the frequency,
linking, and budget constraints~\eqref{eq:frequency}--\eqref{eq:budget} of the
CGN model. They differ only in how passenger routes and capacity constraints are
represented.

Since all alternative models sacrifice some modeling accuracy relative to the CGN, the true quality of any solution they produce can only be assessed by fixing the line plan $(y, z)$ and re-solving the passenger routing problem over the CGN. This yields an exact evaluation of the travel time under that line plan, and we adopt this procedure as our standard evaluation method throughout the paper.

\subsection{PTN-Based Models}
\label{sec:sn}

The simplest alternative routes passengers directly over the PTN, without tracking which line they are traveling on. Passenger routes are paths in the PTN, and capacity constraints are enforced per connection rather than per line, yielding a much more compact model.

Formally, let $\mathcal{R}$ denote the set of all feasible paths in $G^\text{PTN}$, with $\mathcal{R}_k \subseteq \mathcal{R}$ the set of paths for OD pair $k$, and $\mathcal{R}(c) \subseteq \mathcal{R}$ the set of paths containing connection $c$. The travel time of path $r$ is $\tau_r := \sum_{c \in r} \tau_c$. Transfer times are not considered, though dwell times can be incorporated into connection travel times. The resulting formulation is:

\begin{mini!}
{}
{ \sum_{r \in \mathcal{R}} \tau_r x_r + \sum_{k \in \mathcal{OD}} \pi_k u_k }
{}{\texttt{PTN}=}
\addConstraint{ \sum_{r \in \mathcal{R}_k} x_r + u_k }{ \ge v_k }{
    \forall k \in \mathcal{OD}}
\addConstraint{ \sum_{r \in \mathcal{R}(c)} x_r }{ \le
    \sum_{\ell \in \mathcal{L}(c)} \kappa_\ell y_\ell \quad }{
    \forall c \in \mathcal{C}}
\addConstraint{ \eqref{eq:frequency}-\eqref{eq:budget}}{}{}
\addConstraint{ u_k }{ \ge 0 }{ \forall k \in \mathcal{OD}}
\addConstraint{ x_r }{ \ge 0 }{ \forall r \in \mathcal{R}}
\addConstraint{ y_\ell }{ \in \mathbb{Z}_{\ge 0} }{ \forall \ell \in
    \mathcal{L}}
\addConstraint{ z_\ell }{ \in \{0, 1\} }{ \forall \ell \in \mathcal{L}.}
\end{mini!}
The main advantage of \texttt{PTN} is its compactness: the routing graph coincides with the PTN itself, making it the most tractable of the models considered. The key drawback is that transfers are ignored entirely: all routes are assessed as if passengers travel without transferring, and no penalty
is incurred for changing lines. This can lead to solutions that appear attractive in the model but impose significant transfer burden on passengers in practice. 

\subsubsection{Direct Route Constraints}
\label{subsec:ldc}
\citet{borndorfer2015metric} propose a way to partially address this limitation by distinguishing between \emph{direct} and \emph{indirect} routes. A route $r \in \mathcal{R}$ is called a \emph{direct connection route} if it can be traversed by a single line without any transfer. We
denote the set of direct connection routes for OD pair $k$ by
$\mathcal{R}^0_k$, and let $\mathcal{R}^0 = \bigcup_{k \in \mathcal{OD}} \mathcal{R}^0_k$. We also let $\mathcal{L}_r$ denote the set of lines supporting direct route $r \in \mathcal{R}^0$.

Passenger routing is then modeled using two types of variables: $x_{r,0} \ge 0$ for flow on direct route $r \in \mathcal{R}^0$, and $x_{r,1} \ge 0$ for flow on route $r \in \mathcal{R}$ involving at least one transfer. Direct routes incur travel time $\tau_{r,0} := \sum_{c \in r} \tau_c$, while indirect routes incur an additional transfer penalty $\sigma > 0$, giving $\tau_{r,1} := \sigma + \sum_{c \in r} \tau_c$. 

Enforcing exact capacity restrictions on direct routes can be done via Benders decomposition, which we describe later. Here, we present a compact formulation that uses a first-order approximation of direct-travel capacity.

Given a direct route $r \in \mathcal{R}^0$, the total capacity available for direct travel equals $\sum_{\ell \in \mathcal{L}_r} \kappa_\ell y_\ell$. A natural capacity constraint is therefore $x_{r,0} \le \sum_{\ell \in
\mathcal{L}_r} \kappa_\ell y_\ell$ for all $r \in \mathcal{R}^0$. This can be strengthened by observing that direct routes sharing a connection with $r$ and whose supporting line set is a subset of $\mathcal{L}_r$ compete for the same capacity. This gives the \textit{basic direct connection constraints} of
\citet{borndorfer2015metric}:
\[
    \sum_{\tilde{r} \in \mathcal{R}^0(c):\, \mathcal{L}_{\tilde{r}}
    \subseteq \mathcal{L}_r} x_{\tilde{r},0} \le \sum_{\ell \in
    \mathcal{L}_r} \kappa_\ell y_\ell \quad \forall r \in \mathcal{R}^0,\ c \in r.
\]
The number of these constraints can be reduced by aggregating over connections and distinct supporting line-sets. For each connection $c \in \mathcal{C}$, let $\mathbf{L}(c) := \{\mathcal{L}_r : r \in \mathcal{R}^0(c)\}$ denote the collection of distinct line-sets supporting direct routes through $c$. The basic direct connection constraints are then equivalent to
\[
    \sum_{r \in \mathcal{R}^0(c):\, \mathcal{L}_r \subseteq L} x_{r,0}
    \le \sum_{\ell \in L} \kappa_\ell y_\ell \quad \forall c \in
    \mathcal{C},\ L \in \mathbf{L}(c).
\]
We use this more compact form to define the \textit{basic direct connection} model:
\begin{mini!}
{}
{ \sum_{r \in \mathcal{R}^0} \tau_{r,0} x_{r,0} + \sum_{r \in \mathcal{R}}
    \tau_{r,1} x_{r,1} + \sum_{k \in \mathcal{OD}} \pi_k u_k }
{}{\texttt{PTN-BDC}=}
\addConstraint{ \sum_{r \in \mathcal{R}^0_k} x_{r,0} + \sum_{r \in
    \mathcal{R}_k} x_{r,1} + u_k }{ \ge v_k }{
    \forall k \in \mathcal{OD}}
\addConstraint{ \sum_{r \in \mathcal{R}^0(c)} x_{r,0} + \sum_{r \in
    \mathcal{R}(c)} x_{r,1} }{ \le \sum_{\ell \in \mathcal{L}(c)}
    \kappa_\ell y_\ell \quad }{ \forall c \in \mathcal{C}}
\addConstraint{ \sum_{r \in \mathcal{R}^0(c):\, \mathcal{L}_r \subseteq L}
    x_{r,0} }{ \le \sum_{\ell \in L} \kappa_\ell y_\ell \quad }{
    \forall c \in \mathcal{C},\ L \in \mathbf{L}(c)
    \label{eq:direct}}
\addConstraint{ \eqref{eq:frequency}-\eqref{eq:budget}}{}{}
\addConstraint{ u_k }{ \ge 0 }{ \forall k \in \mathcal{OD}}
\addConstraint{ x_{r,0}, x_{r,1} }{ \ge 0 }{ \forall r \in \mathcal{R}}
\addConstraint{ y_\ell }{ \in \mathbb{Z}_{\ge 0} }{ \forall \ell \in
    \mathcal{L}}
\addConstraint{ z_\ell }{ \in \{0, 1\} }{ \forall \ell \in \mathcal{L}.}
\end{mini!}

This model does not fully distinguish between direct and indirect passengers. To illustrate, consider the line pool of Figure~\ref{fig:karbFail}, consisting of four lines: $\ell_1 = A\to B\to C\to D$, $\ell_2 = A\to B\to C$, $\ell_3
= B\to C\to D$, and $\ell_4 = B\to C$. Consider one unit of demand each for OD pairs $A\to C$ and $B\to D$. The supporting line sets are $\mathcal{L}_{A\to C} = \{\ell_1, \ell_2\}$ and $\mathcal{L}_{B\to D} = \{\ell_1, \ell_3\}$; neither set is a subset of the other. The basic direct
connection constraints on connection $B\to C$ therefore read
\begin{align*}
    x_{A\to C,\,0} &\le \kappa(y_{\ell_1} + y_{\ell_2}), \\
    x_{B\to D,\,0} &\le \kappa(y_{\ell_1} + y_{\ell_3}).
\end{align*}
Now suppose $\kappa = 1$ and the model selects $y_{\ell_1} = y_{\ell_4} = 1$ and $y_{\ell_2} = y_{\ell_3} = 0$. Both constraints evaluate to $x_{\cdot,\,0} \le 1$, so the model assigns both passengers as direct and no transfer penalties are incurred. However, $\ell_1$ is the only operated line supporting either direct route: $\ell_4$ covers connection $B\to C$ but does not extend to $A$ or $D$ and cannot support $A\to C$ or $B\to D$ directly. Since $\ell_1$ has capacity $\kappa y_{\ell_1} = 1$, at most one passenger can travel directly. The gap arises because $\mathcal{L}_{A\to C}$ and $\mathcal{L}_{B\to D}$ are incomparable, so no single basic direct connection constraint bounds their combined flow on $\ell_1$ over connection $B\to C$.
\definecolor{tol_green}{HTML}{228833}
\definecolor{tol_orange}{HTML}{CCBB44}
\begin{figure}
\centering
\begin{tikzpicture}[
  stn/.style    = {circle, draw, thick, minimum size=8mm, font=\small},
  conn/.style   = {thick},
  line1/.style  = {very thick, tol_blue},
  line2/.style  = {very thick, tol_red},
  line3/.style  = {very thick, tol_green},
  line4/.style  = {very thick, tol_orange},
]
  \node[stn] (A) at (-3,  0) {A};
  \node[stn] (B) at (-1,  0) {B};
  \node[stn] (C) at ( 1,  0) {C};
  \node[stn] (D) at ( 3,  0) {D};

  \draw[conn] (A) -- (B);
  \draw[conn] (B) -- (C);
  \draw[conn] (C) -- (D);

  % Line 2: A->B->C
  \draw[line1, -{Stealth}]
    (-3, 0.9) -- (-1, 0.9) -- (1, 0.9);
  \node[tol_blue, font=\small, above] at (-1, 0.9) {$\ell_2$};

  % Line 3: B->C->D
  \draw[line2, -{Stealth}]
    (-1,-0.9) -- (1,-0.9) -- (3,-0.9);
  \node[tol_red, font=\small, below] at (1,-0.9) {$\ell_3$};

  % Line 1: A->B->C->D
  \draw[line3, -{Stealth}]
    (-3, 1.7) -- (-1, 1.7) -- (1, 1.7) -- (3, 1.7);
  \node[tol_green, font=\small, above] at (0, 1.7) {$\ell_1$};

  % Line 4: B->C
  \draw[line4, -{Stealth}]
    (-1,-1.7) -- (1,-1.7);
  \node[tol_orange, font=\small, below] at (0,-1.7) {$\ell_4$};

\end{tikzpicture}
\caption{Line pool where the basic direct connection model overestimates direct route capacity.}
\label{fig:karbFail}
\end{figure}
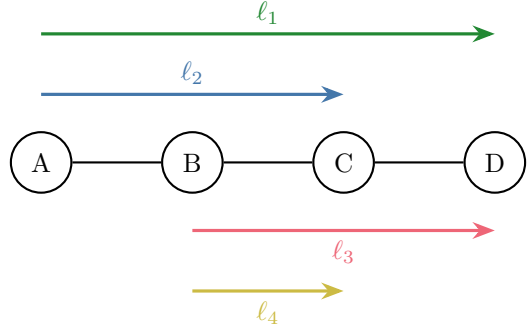

A more accurate treatment can be obtained by extending the static direct
connection constraints~\eqref{eq:direct} with a Benders decomposition, in which the subproblem checks whether the direct passenger flows prescribed by the master problem can be realized by disaggregating them over the lines in the current solution $(y, z)$, generating a cutting plane whenever they cannot. We describe this approach in Section~\ref{sec:sol} and refer to the resulting model as the
\textit{full direct connection} model (\texttt{PTN-FDC}). Even \texttt{PTN-FDC} remains an approximation of the CGN: passengers assigned to indirect routes are charged a single transfer penalty $\sigma$ regardless of how many transfers their route actually requires.

\subsection{Direction-Expanded PTN}
\label{sec:desn}

An alternative approach to partially capturing transfers is to extend the PTN with directional information, offering a middle ground between the CGN and the PTN. Like the CGN, this approach models transfers explicitly; unlike the CGN, it does not track individual lines, resulting in a more compact network. \citet{yao2024railway} apply this idea in what they call the ``extended-direct service network.'' We describe a slight variant with fewer running arcs, which we call the \textit{direction-expanded PTN} (DPTN).

Formally, the DPTN is a directed graph $G^\text{DPTN} = (\mathcal{W}, \mc{E})$. Let $\delta(s) \subseteq \mathcal{C}$ denote the connections incident to station $s$ in the PTN. The node set $\mathcal{W}$ contains one \textit{OD-node} for each station $s \in \mathcal{S}$, and one \textit{travel node} for each station $s \in \mathcal{S}$, incident connection $c \in \delta(s)$, and type $\texttt{t} \in \{\texttt{arr},\texttt{dep}\}$. The arc set $\mc{E}$ consists of four types:
\begin{itemize}
    \item \textit{Running arcs} $\mc{E}^R$: an arc from $(s, c,
        \texttt{dep})$ to $(s', c, \texttt{arr})$ for each connection $c =
        \{s, s'\} \in \mathcal{C}$, representing travel along a connection.
    \item \textit{Dwelling arcs} $\mc{E}^D$: an arc from $(s, c,
        \texttt{arr})$ to $(s, c', \texttt{dep})$ for each pair of
        consecutive connections $(c, c')$ that a line may traverse through
        station $s$, representing continued travel without a transfer.
    \item \textit{Transfer arcs} $\mc{E}^T$: an arc from $(s, c,
        \texttt{arr})$ to $(s, c', \texttt{dep})$ for each pair of connections $c \neq c'$ at station $s$, representing a change of direction at $s$. Infeasible or unlikely transfers may be excluded.
    \item \textit{Boarding and alighting arcs} $\mc{E}^A$: an arc from the OD-node of $s$ to $(s, c, \texttt{dep})$ for each $c \in \delta(s)$, representing boarding, and an arc from $(s, c, \texttt{arr})$ to the OD-node of $s$ for each $c \in \delta(s)$, representing alighting.
\end{itemize}
% \usepackage{tikz}
% \usetikzlibrary{arrows.meta, bending}
% \usepackage{xcolor}
\begin{figure}[htbp]
\centering
\begin{tikzpicture}[
  od/.style    = {circle, draw, thick, minimum size=8mm, font=\small},
  trv/.style   = {rectangle, draw, tol_green, thick, minimum width=1.5cm,
                  minimum height=6mm, font=\scriptsize, rounded corners=2pt},
  trvs/.style  = {rectangle, draw, thick, minimum width=1.5cm,
                  minimum height=6mm, font=\scriptsize, rounded corners=2pt},
  run/.style   = {-{Stealth}, thick, tol_green},
  dwell/.style = {-{Stealth}, thick, tol_green, dotted},
  trf/.style   = {-{Stealth[length=4pt]}, gray, thin, densely dashed},
  board/.style = {-{Stealth[length=3pt]}, lightgray, thin},
]
  % ------ Travel nodes ------
  \node[trv] (AdepAC) at (-4.6,  2.9) {\scriptsize$A\!-\!AC\!-\!\texttt{dep}$};
  \node[trv] (CarrAC) at (-4.0,  1.2) {\scriptsize$C\!-\!AC\!-\!\texttt{arr}$};
  \node[trv] (DdepDE) at ( 4.0,  1.2) {\scriptsize$D\!-\!DE\!-\!\texttt{dep}$};
  \node[trv] (EarrDE) at ( 4.6,  2.9) {\scriptsize$E\!-\!DE\!-\!\texttt{arr}$};
  \node[trv] (BdepBC) at (-4.6, -2.9) {\scriptsize$B\!-\!BC\!-\!\texttt{dep}$};
  \node[trv] (CarrBC) at (-4.0, -1.2) {\scriptsize$C\!-\!BC\!-\!\texttt{arr}$};
  \node[trv] (DdepDF) at ( 4.0, -1.2) {\scriptsize$D\!-\!DF\!-\!\texttt{dep}$};
  \node[trv] (FarrDF) at ( 4.6, -2.9) {\scriptsize$F\!-\!DF\!-\!\texttt{arr}$};
  \node[trvs] (CdepCD) at (-1.5,  0.0) {\scriptsize$C\!-\!CD\!-\!\texttt{dep}$};
  \node[trvs] (DarrCD) at ( 1.5,  0.0) {\scriptsize$D\!-\!CD\!-\!\texttt{arr}$};
  % ------ OD nodes ------
  \node[od] (A) at (-6.0,  2)  {A};
  \node[od] (B) at (-6.0, -2)  {B};
  \node[od] (C) at (-3.5,  0)  {C};
  \node[od] (D) at ( 3.5,  0)  {D};
  \node[od] (E) at ( 6.0,  2)  {E};
  \node[od] (F) at ( 6.0, -2)  {F};
  % ------ Running arcs ------
  \draw[run] (AdepAC) -- (CarrAC);
  \draw[run] (BdepBC) -- (CarrBC);
  \draw[run] (CdepCD) -- (DarrCD);
  \draw[run] (DdepDE) -- (EarrDE);
  \draw[run] (DdepDF) -- (FarrDF);
  % ------ Dwelling arcs ------
  \draw[dwell] (CarrAC) to[bend left=20]  (CdepCD);
  \draw[dwell] (CarrBC) to[bend right=20] (CdepCD);
  \draw[dwell] (DarrCD) to[bend left=20]  (DdepDE);
  \draw[dwell] (DarrCD) to[bend right=20] (DdepDF);
  % ------ Transfer arcs ------
  \draw[trf] (CarrAC) to[bend right=12] (CdepCD);
  \draw[trf] (CarrBC) to[bend left=12]  (CdepCD);
  \draw[trf] (DarrCD) to[bend right=12] (DdepDE);
  \draw[trf] (DarrCD) to[bend left=12]  (DdepDF);
  % ------ Boarding / alighting arcs ------
  \draw[board] (A)      -- (AdepAC);
  \draw[board] (B)      -- (BdepBC);
  \draw[board] (CarrAC) to (C);
  \draw[board] (CarrBC) to (C);
  \draw[board] (C)      to (CdepCD);
  \draw[board] (DarrCD) to (D);
  \draw[board] (D)      to (DdepDE);
  \draw[board] (D)      to (DdepDF);
  \draw[board] (EarrDE) -- (E);
  \draw[board] (FarrDF) -- (F);
  % ------ Legend (two rows, centered around x=0) ------
  \draw[run,   bend left=0] (-3.8,-4.2) -- ++(1.2,0)
    node[right, black, font=\scriptsize] {Running arc};
  \draw[dwell, bend left=0] ( 0.6,-4.2) -- ++(1.2,0)
    node[right, black, font=\scriptsize] {Dwelling arc};
  \draw[trf,   bend left=0] (-3.8,-4.8) -- ++(1.2,0)
    node[right, black, font=\scriptsize] {Transfer arc};
  \draw[board, bend left=0] ( 0.6,-4.8) -- ++(1.2,0)
    node[right, black, font=\scriptsize] {Boarding/alighting arc};
\end{tikzpicture}
\caption{The DPTN for the example network.
  Circles are OD-nodes, rectangles are travel nodes
  $(s\text{-}c\text{-}\texttt{dep/arr})$.}
\label{fig:dptn}
\end{figure}
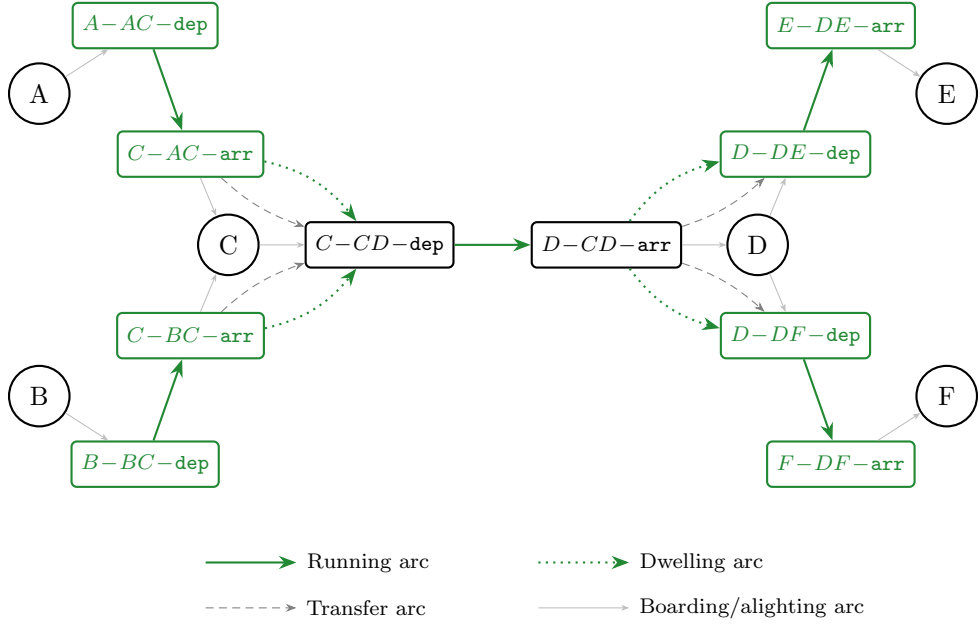
Figure~\ref{fig:dptn} illustrates the DPTN for the PTN of
Figure~\ref{fig:ptn}, restricted to travel from left to right. At stations $C$ and $D$, arrival and departure nodes are connected by two parallel arcs: a dwelling arc, representing a passenger continuing on the same line, and a transfer arc, representing a change of line. Naturally, a dwelling arc may only be used if a line is operated that
traverses both of its incident connections.

We let $\mc{E}^{\text{RD}} = \mc{E}^R \cup \mc{E}^D$ denote
the set of running and dwelling arcs. Each arc $b \in \mc{E}$ has an associated travel time $\tau_b$, and for each $b \in \mc{E}^{\text{RD}}$ we let $\mathcal{L}(b) \subseteq \mathcal{L}$ denote the set of lines that traverse the running or dwelling arc corresponding to $b$. Letting $\mathcal{Q}$ denote the set of feasible passenger paths in the $\text{DPTN}$, with $\mathcal{Q}_k \subseteq \mathcal{Q}$ the set of paths for OD pair $k$ and $\mathcal{Q}(b) \subseteq
\mathcal{Q}$ the set of paths containing arc $b$, the DPTN formulation is:
\begin{mini!}
{}
{ \sum_{q \in \mathcal{Q}} \tau_q x_q + \sum_{k \in \mathcal{OD}} \pi_k u_k }
{}{\texttt{DPTN}=}
\addConstraint{ \sum_{q \in \mathcal{Q}_k} x_q + u_k }{ \ge v_k }{
    \forall k \in \mathcal{OD}}
\addConstraint{ \sum_{q \in \mathcal{Q}(b)} x_q }{ \le
    \sum_{\ell \in \mathcal{L}(b)} \kappa_\ell y_\ell \quad }{
    \forall b \in \mc{E}^{\text{RD}}}
\addConstraint{ \eqref{eq:frequency}-\eqref{eq:budget}}{}{}
\addConstraint{ u_k }{ \ge 0 }{ \forall k \in \mathcal{OD}}
\addConstraint{ x_q }{ \ge 0 }{ \forall q \in \mathcal{Q}}
\addConstraint{ y_\ell }{ \in \mathbb{Z}_{\ge 0} }{ \forall \ell \in
    \mathcal{L}}
\addConstraint{ z_\ell }{ \in \{0, 1\} }{ \forall \ell \in \mathcal{L}.}
\end{mini!}

As illustrated by \citet{yao2024railway}, the DPTN models transfers exactly at a station when the lines involved travel in different directions, or when one of the lines terminates there. However, when both lines traverse the same outgoing connection, the model underestimates the number of required transfers. We illustrate this in more detail in Section~\ref{ss:comp}. 

Capacity constraints are enforced on all running and dwelling arcs, which is somewhat weaker than the per-line capacity constraints of the CGN, since the DPTN does not track which line a passenger is traveling on. Nevertheless,
the DPTN avoids the need to enumerate lines in the routing graph. This results in a network that is larger than the PTN but smaller than the CGN. Crucially, the number of feasible passenger routes is far smaller than in the CGN, since routes are defined over directions rather than over individual lines.

\subsection{Model Hierarchy}
\label{ss:comp}

The models form a natural hierarchy in terms of solution quality. Letting $z(\cdot)$ denote the optimal objective value, it follows immediately that
\[
    z(\texttt{PTN}) \le z(\texttt{PTN-BDC}) \le z(\texttt{PTN-FDC}) \le z(\texttt{CGN})
    \quad \text{and} \quad
    z(\texttt{PTN}) \le z(\texttt{DPTN}) \le z(\texttt{CGN}).
\]
The first chain reflects the progressive tightening of the direct connection relaxation: \texttt{PTN-FDC} dominates \texttt{PTN-BDC} because the Benders cuts it generates are valid inequalities that \texttt{PTN-BDC} may violate.
The direct
connection models \texttt{PTN-BDC} and \texttt{PTN-FDC}, which we jointly refer to
as \texttt{PTN-DC}, are not directly comparable with \texttt{DPTN}: neither
dominates the other in general, as the following examples illustrate.

We stress that this hierarchy holds only for the models as stated, with exhaustive route sets. In practice, the (D)PTN-based solvers enumerate only a restricted set of routes. Since restricting the route set can only increase a model's optimal value, the resulting objective is no longer guaranteed to lower-bound z(\texttt{CGN}): if a model's optimal passenger routing is absent from the enumerated set, the restricted model may overestimate the true cost. Combined with the underestimation already inherent in the approximate models, this means the sign of a solver's estimation error cannot be predicted a priori. We return to this in Section~\ref{sec:detailed}, where we analyze the direction and magnitude of these errors.

\begin{figure}[htbp]
\centering

\begin{subfigure}[c]{0.28\linewidth}
\centering
\scalebox{0.75}{%
\begin{tikzpicture}[
  stn/.style   = {circle, draw, thick, minimum size=8mm, font=\small},
  conn/.style  = {thick},
  line1/.style = {very thick, tol_blue},
  line2/.style = {very thick, tol_red},
]
  % Phantom extents: center the box on y=0 and match panel (b)'s span
  % (diagram +2.9 above, legend down to -4.8 below) so C and D align.
  \path (0,-4.8) (0,4.8);

  \node[stn] (A) at (-2, 1) {A};
  \node[stn] (B) at (-2,-1) {B};
  \node[stn] (C) at (-1, 0) {C};
  \node[stn] (D) at ( 1, 0) {D};
  \node[stn] (E) at ( 2, 1) {E};
  \node[stn] (F) at ( 2,-1) {F};
  \draw[conn] (A) -- (C);
  \draw[conn] (B) -- (C);
  \draw[conn] (C) -- (D);
  \draw[conn] (D) -- (E);
  \draw[conn] (D) -- (F);
  \draw[line1, -{Stealth}]
    (-2, 2.0) -- (-0.9, 0.75) -- (0.9, 0.75) -- (2, 2.0);
  \node[tol_blue, font=\small, above] at (0, 0.75) {$\ell_1$};
  \draw[line2, -{Stealth}]
    (-2,-2) -- (-1,-0.75) -- (1,-0.75) -- (2,-2);
  \node[tol_red, font=\small, below] at (0,-0.75) {$\ell_2$};
\end{tikzpicture}}
\caption{PTN and lines}
\label{fig:ptn-ex1}
\end{subfigure}
\hspace{1cm}
\begin{subfigure}[c]{0.60\linewidth}
\centering
\scalebox{0.75}{%
\begin{tikzpicture}[
  od/.style     = {circle, draw, thick, minimum size=8mm, font=\small},
  trv/.style    = {rectangle, draw, tol_green, thick, minimum width=1.5cm,
                   minimum height=6mm, font=\scriptsize, rounded corners=4pt},
  trvs/.style   = {rectangle, draw, thick, minimum width=1.5cm,
                   minimum height=6mm, font=\scriptsize, rounded corners=4pt},
  run/.style    = {-{Stealth}, thick, tol_green},
  dwell/.style  = {-{Stealth}, thick, tol_green, dotted},
  trf/.style    = {-{Stealth[length=4pt]}, gray, thin, densely dashed},
  board/.style  = {-{Stealth[length=3pt]}, lightgray, thin},
  lbl/.style    = {font=\scriptsize, fill=white, inner sep=1pt},
]
  % Phantom extent: balance the legend (reaching to y=-4.8) with an equal
  % upward extent so the box is centered on y=0, aligning C and D with panel (a).
  \path (0,-4.8) (0,4.8);

  \node[trv]  (AdepAC) at (-4.6,  2.9) {\scriptsize$A\!-\!AC\!-\!\texttt{dep}$};
  \node[trv]  (CarrAC) at (-4.0,  1.2) {\scriptsize$C\!-\!AC\!-\!\texttt{arr}$};
  \node[trv]  (DdepDE) at ( 4.0,  1.2) {\scriptsize$D\!-\!DE\!-\!\texttt{dep}$};
  \node[trv]  (EarrDE) at ( 4.6,  2.9) {\scriptsize$E\!-\!DE\!-\!\texttt{arr}$};
  \node[trv]  (BdepBC) at (-4.6, -2.9) {\scriptsize$B\!-\!BC\!-\!\texttt{dep}$};
  \node[trv]  (CarrBC) at (-4.0, -1.2) {\scriptsize$C\!-\!BC\!-\!\texttt{arr}$};
  \node[trv]  (DdepDF) at ( 4.0, -1.2) {\scriptsize$D\!-\!DF\!-\!\texttt{dep}$};
  \node[trv]  (FarrDF) at ( 4.6, -2.9) {\scriptsize$F\!-\!DF\!-\!\texttt{arr}$};
  \node[trvs] (CdepCD) at (-1.5,  0.0) {\scriptsize$C\!-\!CD\!-\!\texttt{dep}$};
  \node[trvs] (DarrCD) at ( 1.5,  0.0) {\scriptsize$D\!-\!CD\!-\!\texttt{arr}$};
  \node[od] (A) at (-6.0,  2)  {A};
  \node[od] (B) at (-6.0, -2)  {B};
  \node[od] (C) at (-3.5,  0)  {C};
  \node[od] (D) at ( 3.5,  0)  {D};
  \node[od] (E) at ( 6.0,  2)  {E};
  \node[od] (F) at ( 6.0, -2)  {F};
  \draw[run] (AdepAC) -- node[lbl, xshift=-6pt, yshift= 4pt] {1} (CarrAC);
  \draw[run] (BdepBC) -- node[lbl, xshift=-6pt, yshift=-4pt] {1} (CarrBC);
  \draw[run] (CdepCD) -- node[lbl, above]                    {2} (DarrCD);
  \draw[run] (DdepDE) -- node[lbl, xshift= 6pt, yshift= 4pt] {1} (EarrDE);
  \draw[run] (DdepDF) -- node[lbl, xshift= 6pt, yshift=-4pt] {1} (FarrDF);
  \draw[dwell] (CarrAC) to[bend left=20]  node[lbl, above] {1} (CdepCD);
  \draw[dwell] (CarrBC) to[bend right=20] node[lbl, below] {1} (CdepCD);
  \draw[dwell] (DarrCD) to[bend left=20]  node[lbl, above] {1} (DdepDE);
  \draw[dwell] (DarrCD) to[bend right=20] node[lbl, below] {1} (DdepDF);
  \draw[trf] (CarrAC) to[bend right=12] (CdepCD);
  \draw[trf] (CarrBC) to[bend left=12]  (CdepCD);
  \draw[trf] (DarrCD) to[bend right=12] (DdepDE);
  \draw[trf] (DarrCD) to[bend left=12]  (DdepDF);
  \draw[board] (A)      -- (AdepAC);
  \draw[board] (B)      -- (BdepBC);
  \draw[board] (CarrAC) to (C);
  \draw[board] (CarrBC) to (C);
  \draw[board] (C)      to (CdepCD);
  \draw[board] (DarrCD) to (D);
  \draw[board] (D)      to (DdepDE);
  \draw[board] (D)      to (DdepDF);
  \draw[board] (EarrDE) -- (E);
  \draw[board] (FarrDF) -- (F);
  \draw[run,   bend left=0] (-3.8,-4.2) -- ++(1.2,0)
    node[right, black, font=\scriptsize] {Running arc};
  \draw[dwell, bend left=0] ( 0.6,-4.2) -- ++(1.2,0)
    node[right, black, font=\scriptsize] {Dwelling arc};
  \draw[trf,   bend left=0] (-3.8,-4.8) -- ++(1.2,0)
    node[right, black, font=\scriptsize] {Transfer arc};
  \draw[board, bend left=0] ( 0.6,-4.8) -- ++(1.2,0)
    node[right, black, font=\scriptsize] {Boarding/alighting arc};
\end{tikzpicture}}
\caption{DPTN. Arc labels indicate the number of lines traversing each
  running or dwelling arc.}
\label{fig:dptn-ex1b}
\end{subfigure}

\caption{Example~1 with lines $\ell_1 = A\to C\to D\to E$ and
  $\ell_2 = B\to C\to D\to F$.}
\label{fig:dptn-ex1}
\end{figure}
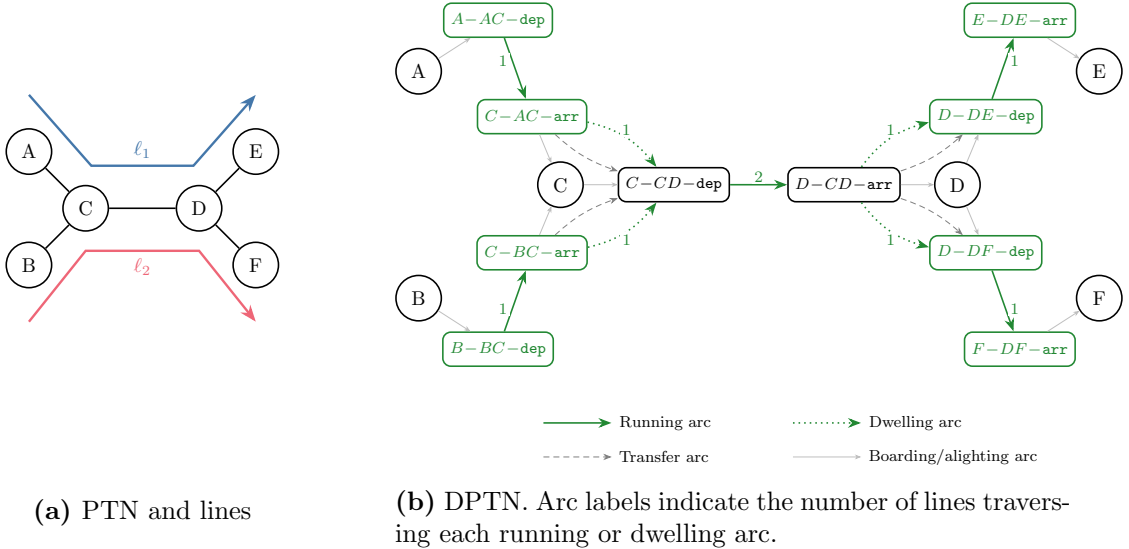

\paragraph{Example 1: \texttt{DPTN} weaker than \texttt{PTN-DC}.}
Consider the example PTN with lines $\ell_1 = A\to C\to D\to E$ and $\ell_2
= B\to C\to D\to F$, shown in Figure~\ref{fig:ptn-ex1}, and four OD pairs $A\to E$, $A\to F$, $B\to E$, and
$B\to F$, each with one passenger. Suppose
the objective is to minimize total transfer penalties, and assume that all
lines have sufficient capacity to carry all passengers. Under the CGN, passengers traveling $A\to F$ and $B\to E$
must transfer, giving $z(\texttt{CGN}) = 2\sigma$. The \texttt{PTN-DC} model
correctly identifies that no direct-route capacity exists for these two OD pairs, so $z(\texttt{PTN-DC}) = 2\sigma$. The DPTN, however, activates all dwelling arcs at $C$ and $D$ (both labeled~1 in Figure~\ref{fig:dptn-ex1b}),
allowing every passenger to travel without transferring, and therefore
underestimates: $z(\texttt{DPTN}) = 0$.
% \usepackage{graphicx} (add to preamble if not present)
% \usepackage{subcaption} (add to preamble if not present)

% =====================================================================
%  Example 1: lines A->C->D->E and B->C->D->F
% =====================================================================
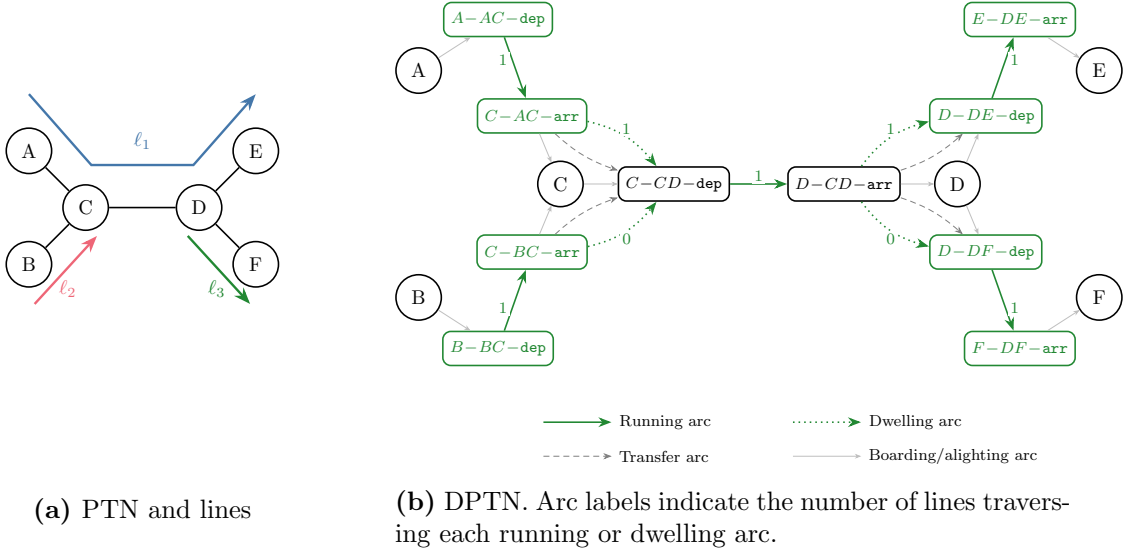
\begin{figure}[htbp]
\centering

\begin{subfigure}[c]{0.28\linewidth}
\centering
\scalebox{0.75}{%
\begin{tikzpicture}[
  stn/.style   = {circle, draw, thick, minimum size=8mm, font=\small},
  conn/.style  = {thick},
  line1/.style = {very thick, tol_blue},
  line2/.style = {very thick, tol_red},
  line3/.style = {very thick, tol_green},
]
  % Phantom extents: center the box on y=0 and match panel (b)'s span
  % (diagram +2.9 above, legend down to -4.8 below) so C and D align.
  \path (0,-4.8) (0,4.8);

  \node[stn] (A) at (-2, 1) {A};
  \node[stn] (B) at (-2,-1) {B};
  \node[stn] (C) at (-1, 0) {C};
  \node[stn] (D) at ( 1, 0) {D};
  \node[stn] (E) at ( 2, 1) {E};
  \node[stn] (F) at ( 2,-1) {F};
  \draw[conn] (A) -- (C);
  \draw[conn] (B) -- (C);
  \draw[conn] (C) -- (D);
  \draw[conn] (D) -- (E);
  \draw[conn] (D) -- (F);
  \draw[line1, -{Stealth}]
    (-2, 2.0) -- (-0.9, 0.75) -- (0.9, 0.75) -- (2, 2.0);
  \node[tol_blue, font=\small, above] at (0, 0.9) {$\ell_1$};
  \draw[line2, -{Stealth}]
    (-1.9,-1.7) -- (-0.8,-0.5);
  \node[tol_red, font=\small, below right] at (-1.6,-1.15) {$\ell_2$};
  \draw[line3, -{Stealth}]
    (0.8,-0.5) -- (1.9,-1.7);
  \node[tol_green, font=\small, below left] at (1.6,-1.15) {$\ell_3$};
\end{tikzpicture}}
\caption{PTN and lines}
\label{fig:ptn-ex2}
\end{subfigure}
\hspace{1cm}
\begin{subfigure}[c]{0.60\linewidth}
\centering
\scalebox{0.75}{%
\begin{tikzpicture}[
  od/.style    = {circle, draw, thick, minimum size=8mm, font=\small},
  trv/.style   = {rectangle, draw, tol_green, thick, minimum width=1.5cm,
                  minimum height=6mm, font=\scriptsize, rounded corners=4pt},
  trvs/.style  = {rectangle, draw, thick, minimum width=1.5cm,
                  minimum height=6mm, font=\scriptsize, rounded corners=4pt},
  run/.style   = {-{Stealth}, thick, tol_green},
  dwell/.style = {-{Stealth}, thick, tol_green, dotted},
  trf/.style   = {-{Stealth[length=4pt]}, gray, thin, densely dashed},
  board/.style = {-{Stealth[length=3pt]}, lightgray, thin},
  lbl/.style   = {font=\scriptsize, fill=white, inner sep=1pt},
]
  % Phantom extent: balance the legend (reaching to y=-4.8) with an equal
  % upward extent so the box is centered on y=0, aligning C and D with panel (a).
  \path (0,-4.8) (0,4.8);

  \node[trv]  (AdepAC) at (-4.6,  2.9) {\scriptsize$A\!-\!AC\!-\!\texttt{dep}$};
  \node[trv]  (CarrAC) at (-4.0,  1.2) {\scriptsize$C\!-\!AC\!-\!\texttt{arr}$};
  \node[trv]  (DdepDE) at ( 4.0,  1.2) {\scriptsize$D\!-\!DE\!-\!\texttt{dep}$};
  \node[trv]  (EarrDE) at ( 4.6,  2.9) {\scriptsize$E\!-\!DE\!-\!\texttt{arr}$};
  \node[trv]  (BdepBC) at (-4.6, -2.9) {\scriptsize$B\!-\!BC\!-\!\texttt{dep}$};
  \node[trv]  (CarrBC) at (-4.0, -1.2) {\scriptsize$C\!-\!BC\!-\!\texttt{arr}$};
  \node[trv]  (DdepDF) at ( 4.0, -1.2) {\scriptsize$D\!-\!DF\!-\!\texttt{dep}$};
  \node[trv]  (FarrDF) at ( 4.6, -2.9) {\scriptsize$F\!-\!DF\!-\!\texttt{arr}$};
  \node[trvs] (CdepCD) at (-1.5,  0.0) {\scriptsize$C\!-\!CD\!-\!\texttt{dep}$};
  \node[trvs] (DarrCD) at ( 1.5,  0.0) {\scriptsize$D\!-\!CD\!-\!\texttt{arr}$};
  \node[od] (A) at (-6.0,  2)  {A};
  \node[od] (B) at (-6.0, -2)  {B};
  \node[od] (C) at (-3.5,  0)  {C};
  \node[od] (D) at ( 3.5,  0)  {D};
  \node[od] (E) at ( 6.0,  2)  {E};
  \node[od] (F) at ( 6.0, -2)  {F};
  \draw[run] (AdepAC) -- node[lbl, xshift=-6pt, yshift= 4pt] {1} (CarrAC);
  \draw[run] (BdepBC) -- node[lbl, xshift=-6pt, yshift=-4pt] {1} (CarrBC);
  \draw[run] (CdepCD) -- node[lbl, above]                    {1} (DarrCD);
  \draw[run] (DdepDE) -- node[lbl, xshift= 6pt, yshift= 4pt] {1} (EarrDE);
  \draw[run] (DdepDF) -- node[lbl, xshift= 6pt, yshift=-4pt] {1} (FarrDF);
  \draw[dwell] (CarrAC) to[bend left=20]  node[lbl, above] {1} (CdepCD);
  \draw[dwell] (CarrBC) to[bend right=20] node[lbl, below] {0} (CdepCD);
  \draw[dwell] (DarrCD) to[bend left=20]  node[lbl, above] {1} (DdepDE);
  \draw[dwell] (DarrCD) to[bend right=20] node[lbl, below] {0} (DdepDF);
  \draw[trf] (CarrAC) to[bend right=12] (CdepCD);
  \draw[trf] (CarrBC) to[bend left=12]  (CdepCD);
  \draw[trf] (DarrCD) to[bend right=12] (DdepDE);
  \draw[trf] (DarrCD) to[bend left=12]  (DdepDF);
  \draw[board] (A)      -- (AdepAC);
  \draw[board] (B)      -- (BdepBC);
  \draw[board] (CarrAC) to (C);
  \draw[board] (CarrBC) to (C);
  \draw[board] (C)      to (CdepCD);
  \draw[board] (DarrCD) to (D);
  \draw[board] (D)      to (DdepDE);
  \draw[board] (D)      to (DdepDF);
  \draw[board] (EarrDE) -- (E);
  \draw[board] (FarrDF) -- (F);
  \draw[run,   bend left=0] (-3.8,-4.2) -- ++(1.2,0)
    node[right, black, font=\scriptsize] {Running arc};
  \draw[dwell, bend left=0] ( 0.6,-4.2) -- ++(1.2,0)
    node[right, black, font=\scriptsize] {Dwelling arc};
  \draw[trf,   bend left=0] (-3.8,-4.8) -- ++(1.2,0)
    node[right, black, font=\scriptsize] {Transfer arc};
  \draw[board, bend left=0] ( 0.6,-4.8) -- ++(1.2,0)
    node[right, black, font=\scriptsize] {Boarding/alighting arc};
\end{tikzpicture}}
\caption{DPTN. Arc labels indicate the number of lines traversing each running or dwelling arc.}
\label{fig:dptn-ex2b}
\end{subfigure}

\caption{Example~2 with lines $\ell_1 = A\to C\to D\to E$,
  $\ell_2 = B\to C$, and $\ell_3 = D\to F$.}
\label{fig:dptn-ex2}
\end{figure}

\paragraph{Example 2: \texttt{PTN-DC} weaker than \texttt{DPTN}.}
Now consider three lines: $\ell_1 = A\to C\to D\to E$, $\ell_2 = B\to C$,
and $\ell_3 = D\to F$, as shown in Figure~\ref{fig:ptn-ex2}. Only
the OD pair $A\to E$ has a direct route; passengers traveling $B\to F$ must
transfer twice. The CGN captures this exactly: $z(\texttt{CGN}) = 4\sigma$.
The DPTN likewise models both transfers correctly, giving $z(\texttt{DPTN}) =
4\sigma$. This is because the dwelling arcs $C\!-\!BC\!-\!\texttt{arr} \to
C\!-\!CD\!-\!\texttt{dep}$ and $D\!-\!CD\!-\!\texttt{arr} \to
D\!-\!DF\!-\!\texttt{dep}$ are both labeled~0 in
Figure~\ref{fig:dptn-ex2b}, forcing passengers to use transfer arcs at
both $C$ and $D$. The \texttt{PTN-DC} model, however, charges at most one
transfer penalty per passenger regardless of the number of transfers made, and
so underestimates: $z(\texttt{PTN-DC}) = 3\sigma$.
% \usepackage{graphicx} (add to preamble if not present)
% \usepackage{subcaption} (add to preamble if not present)
% \usepackage{adjustbox} (add to preamble if not present)

% =====================================================================
%  Example 2: lines A->C->D->E, B->C, and D->F
% =====================================================================

\begin{table}[h]
\centering
\caption{Comparison of network models for line planning.}
\label{tab:models}
\resizebox{\textwidth}{!}{
\begin{tabular}{lllll}
\hline
Model & \centering Transfer penalties & Capacity & Nodes & Arcs \\
\hline
\texttt{PTN}     & \centering None                                     & Per connection & $O(|\mathcal{S}|)$              & $O(|\mathcal{C}|)$              \\
\multirow{2}{*}{\texttt{PTN-DC}} & \centering Approximate direct/indirect split (\texttt{BDC})       & \multirow{2}{*}{Per connection} & \multirow{2}{*}{$O(|\mathcal{S}|)$}              & \multirow{2}{*}{$O(|\mathcal{C}|)$}              \\
& \centering Exact direct/indirect split (\texttt{FDC})              &  &              &             \\
\texttt{DPTN}    & \centering Between terminating/diverging lines & Per connection & $O(|\mathcal{S}|+|\mathcal{C}|)$& $O(|\mathcal{C}|\Delta)$        \\
\texttt{CGN}     & \centering Exact                                    & Per line       & $O(|\mathcal{S}||\mathcal{L}|)$ & $O(|\mathcal{C}||\mathcal{L}|)$ \arraybackslash\\
\hline
\end{tabular}
}
\end{table}

Table~\ref{tab:models} summarizes the key trade-offs between the models,
where $\Delta$ denotes the maximum station degree in the PTN\@. The models \texttt{PTN},
\texttt{PTN-BDC}, and \texttt{PTN-FDC} all share the same routing graph ---
the direct connection constraints and Benders cuts add constraints involving
line variables, but no additional nodes or arcs. They differ in how accurately
they enforce the direct/indirect split: \texttt{PTN-BDC} uses static
constraints that may overestimate the capacity available for direct travel,
while \texttt{PTN-FDC} enforces this exactly via Benders decomposition.
Neither model can charge more than one transfer penalty per passenger,
however, so both underestimate transfer costs when a passenger's route
requires multiple transfers. The DPTN is larger by a factor of $\Delta$,
reflecting the per-direction expansion of each station, but remains
substantially more compact than the CGN when $|\mathcal{L}|$ is large.
Unlike the PTN-based models, it can charge multiple transfer penalties per
passenger, but only when the lines involved terminate or diverge at the
transfer station.
\section{Solution Approaches}
\label{sec:sol}

All presented line planning models share the same structure: a set of integer line frequency variables $y_\ell$ and binary activation variables $z_\ell$, coupled with a (potentially large) set of continuous passenger flow variables. For the PTN-based models, the network is small enough that all passenger routes can be enumerated a priori, and the resulting mixed-integer program is solved directly with a general-purpose MIP solver. Two of the models require additional
machinery. First, for the direct connection model, the direct passenger flows prescribed by the master problem must be verified to be realizable given the current line plan; this is handled via Benders decomposition. Second, for \texttt{CGN}, the number of feasible passenger routes is too large to enumerate explicitly, so routes are generated dynamically via branch-and-price. \texttt{DPTN} occupies a middle ground: its routing graph is larger than the PTN but far smaller than the CGN, so it is not obvious whether a priori route enumeration or branch-and-price is preferable. We therefore implement both and compare them empirically in Section~\ref{sec:experiments}.

\subsection{Benders Decomposition for \texttt{PTN-FDC}}
\label{sec:benders-dr}

\texttt{PTN-BDC} approximates the capacity available for direct travelers using the direct connection constraints~\eqref{eq:direct}. A more accurate treatment can be obtained via Benders decomposition, which we now describe. Although \citet{borndorfer2015metric} do not frame their approach in these terms, it can be naturally interpreted as such.

\paragraph{Master problem.} The Benders master problem is \texttt{PTN-BDC} augmented with all the cuts generated so far. At each iteration, the master problem is solved to obtain a new line plan $(y^*, z^*)$ and direct passenger flows $x^*_{r,0}$, which are then passed to the subproblem. Since the number of routes in the PTN is manageable, all routes are enumerated a priori and the Benders decomposition can be implemented as a callback within a standard MIP solver.

\paragraph{Disaggregation subproblem.} Given a master solution $(y^*,
x^*_{r,0})$, the question of whether the direct passenger flows $x^*_{r,0}$ can actually be realized, i.e., whether each direct traveler can be assigned to a specific supporting line with sufficient capacity, amounts to solving the following LP, which maximizes the number of passengers that
can travel directly:
\begin{maxi!}
{}
{ \sum_{r \in \mathcal{R}^0} \sum_{\ell \in \mathcal{L}_r} w^\ell_{r,0} }
{}{\texttt{DS}(y^*, x^*_{r,0}) = \quad}
\addConstraint{ \sum_{\ell \in \mathcal{L}_r} w^\ell_{r,0} }{ \le x^*_{r,0}
    \quad }{ \forall r \in \mathcal{R}^0 \quad (\alpha_r \ge 0)
    \label{eq:sub-assign}}
\addConstraint{ \sum_{r \in \mathcal{R}^0(c):\, \ell \in \mathcal{L}_r}
    w^\ell_{r,0} }{ \le \kappa_\ell y^*_\ell \quad }{
    \forall \ell \in \mathcal{L},\ c \in \ell \quad
    (\beta^\ell_c \ge 0) \label{eq:sub-capacity}}
\addConstraint{ w^\ell_{r,0} }{ \ge 0 }{
    \forall r \in \mathcal{R}^0,\ \ell \in \mathcal{L}_r,}
\end{maxi!}
where $w^\ell_{r,0}$ represents the flow of direct passengers on route $r$
assigned to line $\ell$, and $\alpha_r$ and $\beta^\ell_c$ are the dual
variables associated with Constraints~\eqref{eq:sub-assign} and
\eqref{eq:sub-capacity}, respectively. Any optimal solution provides a
feasible passenger routing: passengers not routed directly have their
flow $x^*_{r,0} - \sum_{\ell \in \mathcal{L}_r} w^\ell_{r,0}$ transferred
to $x_{r,1}$.

The current master solution is feasible if and only if the optimal value of \texttt{DS} equals $\sum_{r \in \mathcal{R}^0} x^*_{r,0}$, i.e., all prescribed direct flows can be realized. Suppose this is not the case. Then by strong duality the dual optimum is also strictly less than $\sum_{r \in
\mathcal{R}^0} x^*_{r,0}$, so for the optimal dual solution $(\alpha^*,
\beta^*)$:
\[
    \sum_{r \in \mathcal{R}^0} \alpha^*_r x^*_{r,0} +
    \sum_{\ell \in \mathcal{L}} \sum_{c \in \ell} \beta^{\ell*}_c\,
    \kappa_\ell y^*_\ell < \sum_{r \in \mathcal{R}^0} x^*_{r,0}.
\]
This means that the current line plan violates the inequality
\begin{equation}
    \sum_{\ell \in \mathcal{L}} \sum_{c \in \ell} \beta^{\ell*}_c\,
    \kappa_\ell y_\ell \ge \sum_{r \in \mathcal{R}^0}
    (1 - \alpha_r^*)\, x_{r,0},
    \label{eq:dcmetric}
\end{equation}
which is a \textit{dcmetric inequality} in the sense of
\citet{borndorfer2015metric}. This inequality is valid for any feasible
master solution and is added as a Benders cut, tightening the relaxation and ruling out the current infeasible solution.

\subsection{Branch-and-Price for \texttt{CGN} and \texttt{DPTN}}
\label{sec:bnp}

For \texttt{CGN}, the number of feasible passenger routes is too large to enumerate explicitly and branch-and-price is the method of choice. Optionally strengthening the relaxation with valid inequalities (Section~\ref{sec:valid-ineq}) yields the branch-price-and-cut variant. As noted above, it is not clear a priori whether enumeration or branch-and-price is preferable for \texttt{DPTN}, so we also implement a branch-and-price variant, with the pricing problem solved over the DPTN. Apart from the route definition, the branch-and-price and branch-price-and-cut algorithms described for the CGN carry over to the DPTN without modification.

Branch-and-price is a variant of branch-and-bound in which the LP relaxation at each node of the search tree is solved by column generation rather than by a direct LP solver. At each node, column generation iterates between a \textit{restricted master problem} (RMP), which optimizes over a subset of routes $\bar{\mathcal{P}} \subseteq \mathcal{P}$, and a \textit{pricing problem}, which searches for a route with negative reduced cost to add to the RMP. The process terminates when no such route exists, certifying LP optimality.

\paragraph{Restricted master problem.} The RMP is the LP relaxation of
\texttt{CGN} restricted to $\bar{\mathcal{P}}$:
\begin{mini!}
{}
{ \sum_{p \in \bar{\mathcal{P}}} \tau_p x_p +
    \sum_{k \in \mathcal{OD}} \pi_k u_k }
{}{\texttt{RMP}(\bar{\mathcal{P}}) = \quad}
\addConstraint{ \sum_{p \in \bar{\mathcal{P}}_k} x_p + u_k }{ \ge v_k }{
    \forall k \in \mathcal{OD} \label{eq:rmp-demand}}
\addConstraint{ \sum_{p \in \bar{\mathcal{P}}(a)} x_p }{ \le \kappa_\ell
    y_\ell \quad }{ \forall \ell \in \mathcal{L},\ a \in \ell
    \label{eq:rmp-capacity}}
\addConstraint{ \eqref{eq:frequency}-\eqref{eq:budget}}{}{}
\addConstraint{ u_k }{ \ge 0 }{ \forall k \in \mathcal{OD}}
\addConstraint{ x_p }{ \ge 0 }{ \forall p \in \bar{\mathcal{P}}}
\addConstraint{ y_\ell }{ \ge 0 }{ \forall \ell \in \mathcal{L}}
\addConstraint{ z_\ell }{ \ge 0 }{ \forall \ell \in \mathcal{L}.}
\end{mini!}

\paragraph{Pricing problem.} Let $\mu_k \ge 0$ and $\lambda_a \ge 0$ denote the dual variables associated
with Constraints~\eqref{eq:rmp-demand} and \eqref{eq:rmp-capacity},
respectively. The reduced cost of route $p \in \mathcal{P}_k$ is
\begin{equation}
    \bar{\tau}_p = \tau_p - \mu_k + \sum_{a \in p} \lambda_a.
    \label{eq:reduced-cost}
\end{equation}
For each OD pair $k \in \mathcal{OD}$, the
pricing problem seeks a path $p \in \mathcal{P}_k$ minimizing the reduced
cost~\eqref{eq:reduced-cost}, i.e., a shortest path in the CGN from the
origin to the destination of $k$ with arc weights $\tau_a + \lambda_a$.
Since these weights are non-negative, the pricing problem can be solved
efficiently using Dijkstra's algorithm. If the minimum reduced cost is
negative, the corresponding route is added to $\bar{\mathcal{P}}$. Otherwise, column generation terminates.

\paragraph{Branching.}  We apply the following branching strategy. We first branch on the total frequency of a connection. If all connection frequencies are integer, we branch on the activation variable $z_\ell$ of a fractional line, fixing it to zero or one. If all activation variables are integer but some line frequencies remain fractional, we branch directly on $y_\ell$. All branching decisions involve line variables only, and do not affect the structure of the pricing problem.

To select among candidate variables at each node, we use \emph{semi-strong branching}. Among the ten most fractional candidates according to the current branching rule, we evaluate each by solving the restricted master problem with the candidate branching decision applied, but without performing additional column generation iterations. The candidate that leads to the highest weighted increase in objective value across the two child nodes is selected.

\paragraph{Primal heuristics.}
Branching alone typically only yields integer solutions deep in the branch-and-price tree. As such, we use two different primal heuristics to find good feasible solutions early on.

At the root node, after solving the LP relaxation, we apply a diving heuristic to obtain an initial integer solution. At each iteration, we select the line with the highest fractional activation value and round it up, provided the line budget is not exceeded. If no such line exists, we instead round the most fractional line down. Additionally, after each fixing decision, any line whose activation would cause a budget violation is immediately fixed to zero. After each fixing step, the RMP is re-solved using column generation before the next rounding decision is made. Since every iteration either activates a line within the remaining budget or deactivates one, and the line pool is finite, this procedure is guaranteed to terminate with a feasible integer solution.

At every tenth node of the branch-and-price tree, we solve the current RMP as a MIP, restricted to the columns currently included in the RMP and imposing a time limit of ten seconds. Although this heuristic is not guaranteed to find a feasible solution within the time limit, it is inexpensive and frequently yields good integer solutions early in the search.

\subsubsection{Strong Inequalities}
\label{sec:valid-ineq}
The capacity constraints~\eqref{eq:capacity} bound the total passenger flow on each arc of line~$\ell$ by its aggregate capacity $\kappa_\ell y_\ell$, and are known to lead to weak LP relaxations. In particular, they do not prevent a single OD pair $k$ from being assigned more flow on arc $a \in \ell$ than the capacity that line $\ell$ can provide to that pair, namely $v_k y_\ell$. This observation motivates the use of \emph{strong} (or \emph{disaggregated}) capacity inequalities, originally introduced in the network design literature \citep{magnanti1984network} and previously used in line planning by \citet{yao2024railway}. The strong inequalities read
\begin{equation}
    \sum_{p \in \bar{\mathcal{P}}_{k}(a)} x_p  \le v_k y_\ell 
    \quad \forall k \in \mathcal{OD},\ \ell \in \mathcal{L},\ a \in \ell.
    \label{eq:strong-inequality}
\end{equation}
These inequalities dominate the aggregate capacity constraints~\eqref{eq:capacity} for individual OD pairs and can significantly strengthen the LP relaxation. Since the number of possible strong inequalities is large, we do not include all of them ex ante but generate them in a cutting plane fashion.

\paragraph{Separation.} Checking whether inequality~\eqref{eq:strong-inequality} is violated for some  triple $(k, \ell, a)$ requires only a comparison of the 
left- and right-hand sides, making separation straightforward. Since the number of 
triples $(k, \ell, a)$ can be large, we limit the number of cuts added per column 
generation iteration: for each OD pair $k$, we identify the most violated arc $a^*_k$ 
across all lines and add the corresponding cut only if the relative violation exceeds 
$1\%$.

\paragraph{Line activation inequalities.} As with the aggregate capacity constraints, the linking constraints~\eqref{eq:linking} suffer from a weak LP relaxation: they do not prevent individual pairs from being routed over an arc of a line that is not activated. We therefore use the analogous disaggregated inequalities
\begin{equation}
    \sum_{p \in \bar{\mathcal{P}}_{k}(a)} x_p  \le v_k z_\ell
    \quad \forall k \in \mathcal{OD},\ \ell \in \mathcal{L},\ a \in \ell.
    \label{eq:strong-inequality-activation}
\end{equation}
These inequalities can be separated in an identical way as the strong capacity inequalities.

\subsection{Demand Sparsification}
\label{sec:preprocessing}

For large-scale instances, reducing the number of OD pairs through 
preprocessing can significantly speed up the solution process. We employ 
two techniques.

\paragraph{OD splitting.} Given an OD pair $k$ with $v_k$ passengers 
traveling from $o_k$ to $d_k$, we compute the shortest path in the DPTN. 
If this path contains a transfer at some intermediate station $t_k$, we 
remove the original pair and replace it with two sub-pairs $o_k \to t_k$ 
and $t_k \to d_k$, each with volume $v_k$. This process is repeated until 
all OD pairs have a direct shortest path. Beyond reducing the number of OD 
pairs, splitting shortens average route lengths and reduces the number of 
transfers in the modified instance, simplifying the solution process.

\paragraph{Through-station aggregation.} We classify each station in the 
PTN as a \emph{junction station} (degree three or higher) or a \emph{through 
station} (degree one or two). Since passengers originating at or destined 
for neighboring through stations follow similar routes --- first toward a 
bounding junction station, then onward to their destination --- we aggregate 
OD pairs whose origins or destinations are neighboring through stations. 
This aggregation is repeated until at most one through station between any 
two junction stations retains nonzero demand.
\section{Computational Results}
\label{sec:experiments}
This section presents a comprehensive computational study comparing the models and their solver variants. We first describe the benchmark instances, derived from the Dutch and Swiss railway networks, and define a total of 16 line planning solvers in Section~\ref{sec:instances}. We then divide the solvers into four families, and analyze the relative performance at the family-level in Section~\ref{sec:families}, before examining the detailed performance of the individual variants in Section~\ref{sec:detailed}. Finally, we assess the impact of the demand sparsification techniques in Section~\ref{sec:demSpar}. Our code, instances and best solutions are publicly available at \url{https://github.com/rn-van-lieshout/line-planning-benchmark}.

\subsection{Instances}
\label{sec:instances}
The instances are based on the railway networks of the Netherlands and Switzerland, both of which serve over a million passengers per day \citep{ns_reizigersgedrag,sbb_transportation}. The Dutch railway infrastructure was obtained directly from NS, the country's largest railway operator. NS operates two train types on this network: \emph{intercity} lines, which stop only at major stations, and \emph{sprinter} lines, which stop at every station. Beyond the set of stations and tracks, the data specifies which stations may serve as terminal stations for intercity and sprinter lines, which track sequences are infeasible for through service, and at which stations an intercity line may transition into a sprinter line or vice versa. The Swiss infrastructure was obtained from \href{http://data.sbb.ch}{data.sbb.ch}, with through stations that have fewer than three platforms removed.  Major stations were designated as potential terminal stations.

\begin{figure}[h]
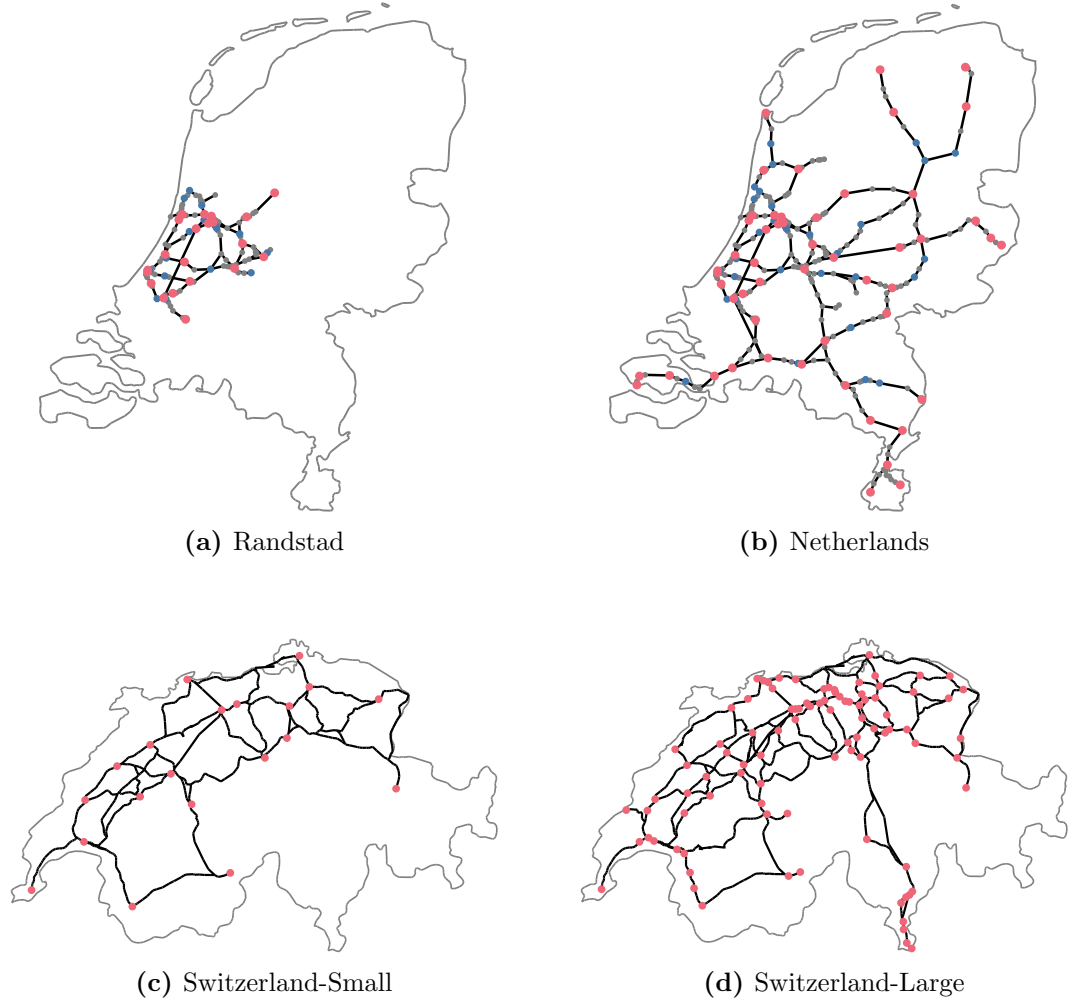

\begin{minipage}[b]{0.5\textwidth}
    \centering
    \resizebox{0.8\textwidth}{!}{\input{figs/NederlandKlein}}
    \subcaption{Randstad}
  \end{minipage}
  \hfill
   \begin{minipage}[b]{0.5\textwidth}
    \centering
    \resizebox{0.8\textwidth}{!}{\input{figs/Nederland}}
    \subcaption{Netherlands}
  \end{minipage}
  \begin{minipage}[b]{0.5\textwidth}
    \centering
    \vspace{1cm}
    \resizebox{0.9\textwidth}{!}{\input{figs/ZwitserlandKlein}}
    \subcaption{Switzerland-Small}
  \end{minipage}
  \hfill
  \begin{minipage}[b]{0.5\textwidth}
    \centering
    \vspace{1cm}
    \resizebox{0.9\textwidth}{!}{\input{figs/Zwitserland}}
    \subcaption{Switzerland-Large}
  \end{minipage}
  \caption{Rail networks in the Netherlands and Switzerland. In the Dutch networks, large red dots correspond to intercity stations, small gray dots to sprinter stations, and the medium-sized blue dots to stations where intercity lines may transition into sprinter lines and vice versa.}
  \label{fig:maps}
\end{figure}

For the Dutch network, we consider four instances: Randstad-IC, Randstad-IC-SPR, Netherlands-IC, and Netherlands-IC-SPR. The Randstad instances cover the densely populated western region of the Netherlands, while the Netherlands instances span the entire national network. The IC variants include only intercity lines, whereas the IC-SPR variants additionally include sprinter lines. For Switzerland, we consider two instances: Switzerland-Small, restricted to the larger cities, and Switzerland-Large, covering the entire network. The networks are visualized in Figure~\ref{fig:maps}. As no demand data was available, we estimated OD matrices using a gravity model based on inter-station distances and city populations.

\textbf{Instance parameters.} For every network, we generate instances by varying five parameters: the line pool (small, medium, or large), the budget level (low, medium, or high), the transfer penalty (5, 15, or 25 minutes), the fixed line cost (low, medium, or high), and the OD matrix (dense or sparse). This yields 162 instances per network, or 972 instances in total.

The line pools differ in how permissive they are with respect to line length. In the small pool, a line's length may deviate by at most 10\% from the corresponding shortest path, and sprinter lines are capped at 60 minutes. The medium and large pools relax these limits to 30\% and 90 minutes, and 50\% and 120 minutes, respectively. The sparse OD matrix is generated from the dense matrix by applying sparsification as in Section~\ref{sec:preprocessing}. Table~\ref{tab:instances} presents for all instances the number of stations and connections, the number of OD pairs in the dense and sparse matrices, and the line pool size for the three variants. 

\begin{table}[h]
\centering
\caption{Instance statistics. $|\mathcal{OD}|$ denotes the number of OD pairs in the dense and sparse matrices. $|\mathcal{L}|$ denotes the line pool size for the small, medium, and large variants.}
\label{tab:instances}
\begin{tabular}{lrrrrrrr}
\toprule
& & & \multicolumn{2}{c}{$|\mathcal{OD}|$} & \multicolumn{3}{c}{$|\mathcal{L}|$} \\
\cmidrule(lr){4-5} \cmidrule(lr){6-8}
Instance & $|\mc{S}|$ & $|\mc{C}|$ & Dense & Sparse & Small & Medium & Large \\
\midrule
Randstad-IC         & 29 &  80 &    416 &    326 &    93 &    140 &    162 \\
Randstad-IC-SPR     & 109 & 324 &  7,534 &  2,969 &   296 &    528 &    662 \\
Netherlands-IC        & 78 & 196 &  2,564 &  2,050 &   590 &  1,266 &  2,364 \\
Netherlands-IC-SPR    & 256 & 750 & 25,510 & 8,965 &   970 &  2,165 &  3,981 \\
\midrule
Switzerland-Small   &  20 &  80 &    380 &    344 &   167 &    252 &    339 \\
Switzerland-Large   &  93 & 302 &  3,418 &  3,002 & 1,201 &  5,312 & 12,359 \\
\bottomrule
\end{tabular}
\end{table}

\textbf{Solver variants.} For models that require a priori route enumeration, we consider three route sets of increasing size, indicated by the suffixes \texttt{-S}, \texttt{-M}, and \texttt{-L}. The small route set (\texttt{-S}) contains, for each OD pair, only those routes whose duration equals that of the shortest path. The medium (\texttt{-M}) and large (\texttt{-L}) sets additionally admit routes with detours of up to 15 and 30 minutes, respectively. For models that rely on column generation, routes are generated dynamically and no such suffix is needed; we instead test variants with and without separation of strong inequalities, referred to as branch-and-price (\texttt{B\&P}) and branch-price-and-cut (\texttt{BP\&C}), respectively. This yields 16 solver variants, which we divide into four families. The \texttt{PTN} family comprises the three enumeration variants \texttt{PTN-S}, \texttt{PTN-M}, and \texttt{PTN-L}. The \texttt{PTN-DC} family collects the direct connection models in both their basic and Benders forms, across all three route sets: \texttt{PTN-BDC-S/M/L} and \texttt{PTN-FDC-S/M/L}. The \texttt{DPTN} family contains the three enumeration variants \texttt{DPTN-S/M/L} together with the two column generation variants \texttt{DPTN-B\&P} and \texttt{DPTN-BP\&C}. Finally, the \texttt{CGN} family consists of the two column generation variants \texttt{CGN-B\&P} and \texttt{CGN-BP\&C}.

\textbf{Implementation details.} We generate the line and route pools by applying a labeling algorithm on the DPTN, which directly encodes infrastructural constraints on train movements: dwelling arcs are included only where the infrastructure permits a train to enter a station via one connection and exit via another. Feasible lines and direct routes then correspond precisely to paths that use only running and dwelling arcs, while paths containing one or more transfer arcs are indirect routes and incur transfer penalties. For the PTN-based models, these DPTN routes are converted to PTN routes and duplicates are removed.

\textbf{Set-up.} We run each of the 16 solver variants on all 972 instances using \texttt{CPLEX22.1.0} with six threads and a one-hour time limit per instance. All column generation algorithms make use of dual smoothing and column management. The pricing problems are solved in parallel on up to eight threads, and semi-strong branching is conducted for at most ten candidates in parallel on at most four threads. All experiments are conducted on cluster nodes with 16GB RAM and an AMD Rome 7H12 processor. 

We always evaluate solution quality of a line plan using the CGN. Since we are interested in analyzing relative solver performance, the upcoming sections mainly present performance profiles, which report for each performance ratio $\rho \geq 1$ the fraction of instances on which the solver finds a feasible solution whose objective lies within a factor~$\rho$ of the best solution found by any solver variant.

\subsection{Comparing Solver Families}
\label{sec:families}

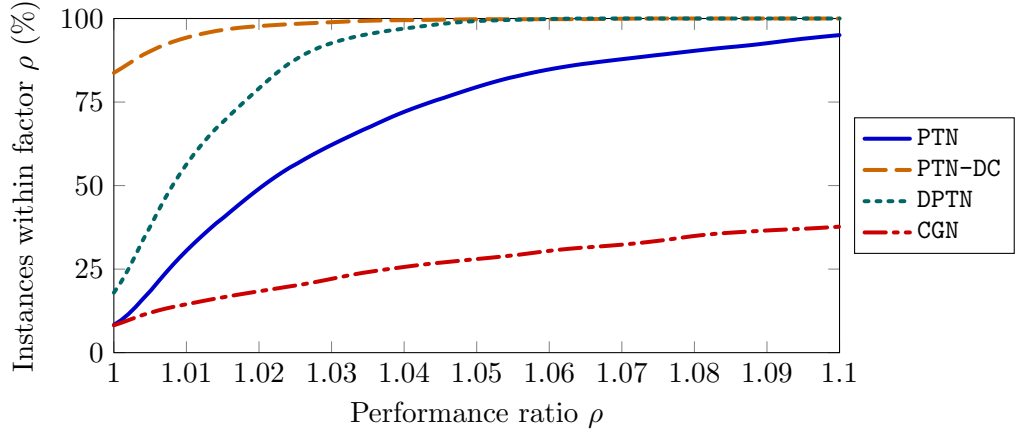
\begin{figure}[t]
\centering
\begin{tikzpicture}
\begin{axis}[
    width=0.75\textwidth,
    height=6cm,
    xmin=1, xmax=1.1000,
    ymin=0, ymax=1,
    ytick={0,0.25,0.5,0.75,1},
    yticklabels={0,25,50,75,100},
    xlabel={Performance ratio $\rho$},
    ylabel={Instances within factor $\rho$ (\%)},
    legend style={at={(1.02,0.5)}, anchor=west, legend columns=1, font=\small},
    legend cell align={left},
    clip=false,
]
\addplot[color=blue!80!black, solid, line width=1.5pt, line cap=round, line join=round] coordinates {
    (1.0000,0.0829)
    (1.0010,0.0971)
    (1.0020,0.1152)
    (1.0030,0.1367)
    (1.0040,0.1607)
    (1.0051,0.1862)
    (1.0061,0.2122)
    (1.0071,0.2376)
    (1.0081,0.2621)
    (1.0091,0.2853)
    (1.0101,0.3073)
    (1.0111,0.3283)
    (1.0121,0.3485)
    (1.0131,0.3681)
    (1.0141,0.3871)
    (1.0152,0.4058)
    (1.0162,0.4242)
    (1.0172,0.4422)
    (1.0182,0.4599)
    (1.0192,0.4772)
    (1.0202,0.4939)
    (1.0212,0.5099)
    (1.0222,0.5251)
    (1.0232,0.5396)
    (1.0242,0.5533)
    (1.0253,0.5665)
    (1.0263,0.5791)
    (1.0273,0.5912)
    (1.0283,0.6028)
    (1.0293,0.6141)
    (1.0303,0.6249)
    (1.0313,0.6355)
    (1.0323,0.6458)
    (1.0333,0.6561)
    (1.0343,0.6664)
    (1.0354,0.6767)
    (1.0364,0.6869)
    (1.0374,0.6968)
    (1.0384,0.7064)
    (1.0394,0.7155)
    (1.0404,0.7241)
    (1.0414,0.7323)
    (1.0424,0.7401)
    (1.0434,0.7478)
    (1.0444,0.7553)
    (1.0455,0.7627)
    (1.0465,0.7701)
    (1.0475,0.7773)
    (1.0485,0.7844)
    (1.0495,0.7914)
    (1.0505,0.7981)
    (1.0515,0.8046)
    (1.0525,0.8108)
    (1.0535,0.8167)
    (1.0545,0.8223)
    (1.0556,0.8276)
    (1.0566,0.8326)
    (1.0576,0.8374)
    (1.0586,0.8419)
    (1.0596,0.8461)
    (1.0606,0.8501)
    (1.0616,0.8538)
    (1.0626,0.8573)
    (1.0636,0.8607)
    (1.0646,0.8638)
    (1.0657,0.8667)
    (1.0667,0.8696)
    (1.0677,0.8723)
    (1.0687,0.8749)
    (1.0697,0.8775)
    (1.0707,0.8800)
    (1.0717,0.8825)
    (1.0727,0.8850)
    (1.0737,0.8876)
    (1.0747,0.8902)
    (1.0758,0.8927)
    (1.0768,0.8953)
    (1.0778,0.8978)
    (1.0788,0.9003)
    (1.0798,0.9028)
    (1.0808,0.9052)
    (1.0818,0.9076)
    (1.0828,0.9099)
    (1.0838,0.9121)
    (1.0848,0.9142)
    (1.0859,0.9163)
    (1.0869,0.9185)
    (1.0879,0.9207)
    (1.0889,0.9232)
    (1.0899,0.9258)
    (1.0909,0.9285)
    (1.0919,0.9313)
    (1.0929,0.9341)
    (1.0939,0.9367)
    (1.0949,0.9393)
    (1.0960,0.9417)
    (1.0970,0.9440)
    (1.0980,0.9462)
    (1.0990,0.9483)
    (1.1000,0.9504)
};
\addlegendentry{\texttt{PTN}}
\addplot[color=orange!80!black, dash pattern=on 8pt off 4pt, line width=1.5pt, line cap=round, line join=round] coordinates {
    (1.0000,0.8374)
    (1.0010,0.8504)
    (1.0020,0.8643)
    (1.0030,0.8781)
    (1.0040,0.8912)
    (1.0051,0.9029)
    (1.0061,0.9133)
    (1.0071,0.9223)
    (1.0081,0.9302)
    (1.0091,0.9371)
    (1.0101,0.9434)
    (1.0111,0.9490)
    (1.0121,0.9540)
    (1.0131,0.9586)
    (1.0141,0.9628)
    (1.0152,0.9664)
    (1.0162,0.9695)
    (1.0172,0.9721)
    (1.0182,0.9742)
    (1.0192,0.9760)
    (1.0202,0.9776)
    (1.0212,0.9791)
    (1.0222,0.9804)
    (1.0232,0.9818)
    (1.0242,0.9830)
    (1.0253,0.9842)
    (1.0263,0.9853)
    (1.0273,0.9863)
    (1.0283,0.9873)
    (1.0293,0.9883)
    (1.0303,0.9893)
    (1.0313,0.9902)
    (1.0323,0.9912)
    (1.0333,0.9920)
    (1.0343,0.9928)
    (1.0354,0.9935)
    (1.0364,0.9940)
    (1.0374,0.9944)
    (1.0384,0.9947)
    (1.0394,0.9949)
    (1.0404,0.9952)
    (1.0414,0.9954)
    (1.0424,0.9957)
    (1.0434,0.9960)
    (1.0444,0.9963)
    (1.0455,0.9966)
    (1.0465,0.9969)
    (1.0475,0.9972)
    (1.0485,0.9975)
    (1.0495,0.9977)
    (1.0505,0.9978)
    (1.0515,0.9979)
    (1.0525,0.9979)
    (1.0535,0.9979)
    (1.0545,0.9979)
    (1.0556,0.9979)
    (1.0566,0.9979)
    (1.0576,0.9979)
    (1.0586,0.9979)
    (1.0596,0.9979)
    (1.0606,0.9980)
    (1.0616,0.9980)
    (1.0626,0.9980)
    (1.0636,0.9980)
    (1.0646,0.9981)
    (1.0657,0.9983)
    (1.0667,0.9985)
    (1.0677,0.9987)
    (1.0687,0.9990)
    (1.0697,0.9993)
    (1.0707,0.9995)
    (1.0717,0.9997)
    (1.0727,0.9998)
    (1.0737,0.9999)
    (1.0747,1.0000)
    (1.0758,1.0000)
    (1.0768,1.0000)
    (1.0778,1.0000)
    (1.0788,1.0000)
    (1.0798,1.0000)
    (1.0808,1.0000)
    (1.0818,1.0000)
    (1.0828,1.0000)
    (1.0838,1.0000)
    (1.0848,1.0000)
    (1.0859,1.0000)
    (1.0869,1.0000)
    (1.0879,1.0000)
    (1.0889,1.0000)
    (1.0899,1.0000)
    (1.0909,1.0000)
    (1.0919,1.0000)
    (1.0929,1.0000)
    (1.0939,1.0000)
    (1.0949,1.0000)
    (1.0960,1.0000)
    (1.0970,1.0000)
    (1.0980,1.0000)
    (1.0990,1.0000)
    (1.1000,1.0000)
};
\addlegendentry{\texttt{PTN-DC}}
\addplot[color=teal!80!black, dash pattern=on 1.5pt off 3pt, line width=1.5pt, line cap=round, line join=round] coordinates {
    (1.0000,0.1796)
    (1.0010,0.2116)
    (1.0020,0.2498)
    (1.0030,0.2923)
    (1.0040,0.3367)
    (1.0051,0.3810)
    (1.0061,0.4236)
    (1.0071,0.4636)
    (1.0081,0.5007)
    (1.0091,0.5349)
    (1.0101,0.5665)
    (1.0111,0.5958)
    (1.0121,0.6228)
    (1.0131,0.6480)
    (1.0141,0.6715)
    (1.0152,0.6938)
    (1.0162,0.7150)
    (1.0172,0.7354)
    (1.0182,0.7555)
    (1.0192,0.7753)
    (1.0202,0.7949)
    (1.0212,0.8142)
    (1.0222,0.8328)
    (1.0232,0.8503)
    (1.0242,0.8663)
    (1.0253,0.8807)
    (1.0263,0.8932)
    (1.0273,0.9040)
    (1.0283,0.9133)
    (1.0293,0.9214)
    (1.0303,0.9285)
    (1.0313,0.9348)
    (1.0323,0.9405)
    (1.0333,0.9457)
    (1.0343,0.9503)
    (1.0354,0.9546)
    (1.0364,0.9585)
    (1.0374,0.9620)
    (1.0384,0.9652)
    (1.0394,0.9682)
    (1.0404,0.9710)
    (1.0414,0.9738)
    (1.0424,0.9766)
    (1.0434,0.9793)
    (1.0444,0.9820)
    (1.0455,0.9846)
    (1.0465,0.9868)
    (1.0475,0.9888)
    (1.0485,0.9904)
    (1.0495,0.9918)
    (1.0505,0.9928)
    (1.0515,0.9937)
    (1.0525,0.9944)
    (1.0535,0.9951)
    (1.0545,0.9957)
    (1.0556,0.9963)
    (1.0566,0.9969)
    (1.0576,0.9974)
    (1.0586,0.9980)
    (1.0596,0.9984)
    (1.0606,0.9989)
    (1.0616,0.9992)
    (1.0626,0.9995)
    (1.0636,0.9997)
    (1.0646,0.9998)
    (1.0657,0.9999)
    (1.0667,1.0000)
    (1.0677,1.0000)
    (1.0687,1.0000)
    (1.0697,1.0000)
    (1.0707,1.0000)
    (1.0717,1.0000)
    (1.0727,1.0000)
    (1.0737,1.0000)
    (1.0747,1.0000)
    (1.0758,1.0000)
    (1.0768,1.0000)
    (1.0778,1.0000)
    (1.0788,1.0000)
    (1.0798,1.0000)
    (1.0808,1.0000)
    (1.0818,1.0000)
    (1.0828,1.0000)
    (1.0838,1.0000)
    (1.0848,1.0000)
    (1.0859,1.0000)
    (1.0869,1.0000)
    (1.0879,1.0000)
    (1.0889,1.0000)
    (1.0899,1.0000)
    (1.0909,1.0000)
    (1.0919,1.0000)
    (1.0929,1.0000)
    (1.0939,1.0000)
    (1.0949,1.0000)
    (1.0960,1.0000)
    (1.0970,1.0000)
    (1.0980,1.0000)
    (1.0990,1.0000)
    (1.1000,1.0000)
};
\addlegendentry{\texttt{DPTN}}
\addplot[color=red!80!black, dash pattern=on 8pt off 3pt on 1.5pt off 3pt, line width=1.5pt, line cap=round, line join=round] coordinates {
    (1.0000,0.0820)
    (1.0010,0.0892)
    (1.0020,0.0971)
    (1.0030,0.1052)
    (1.0040,0.1129)
    (1.0051,0.1199)
    (1.0061,0.1262)
    (1.0071,0.1317)
    (1.0081,0.1365)
    (1.0091,0.1410)
    (1.0101,0.1453)
    (1.0111,0.1495)
    (1.0121,0.1538)
    (1.0131,0.1580)
    (1.0141,0.1621)
    (1.0152,0.1662)
    (1.0162,0.1701)
    (1.0172,0.1738)
    (1.0182,0.1774)
    (1.0192,0.1809)
    (1.0202,0.1843)
    (1.0212,0.1878)
    (1.0222,0.1912)
    (1.0232,0.1947)
    (1.0242,0.1982)
    (1.0253,0.2018)
    (1.0263,0.2055)
    (1.0273,0.2094)
    (1.0283,0.2135)
    (1.0293,0.2178)
    (1.0303,0.2222)
    (1.0313,0.2266)
    (1.0323,0.2309)
    (1.0333,0.2350)
    (1.0343,0.2388)
    (1.0354,0.2424)
    (1.0364,0.2457)
    (1.0374,0.2488)
    (1.0384,0.2518)
    (1.0394,0.2547)
    (1.0404,0.2576)
    (1.0414,0.2604)
    (1.0424,0.2631)
    (1.0434,0.2656)
    (1.0444,0.2680)
    (1.0455,0.2703)
    (1.0465,0.2725)
    (1.0475,0.2747)
    (1.0485,0.2769)
    (1.0495,0.2791)
    (1.0505,0.2812)
    (1.0515,0.2832)
    (1.0525,0.2854)
    (1.0535,0.2876)
    (1.0545,0.2900)
    (1.0556,0.2925)
    (1.0566,0.2953)
    (1.0576,0.2981)
    (1.0586,0.3009)
    (1.0596,0.3036)
    (1.0606,0.3062)
    (1.0616,0.3085)
    (1.0626,0.3106)
    (1.0636,0.3125)
    (1.0646,0.3143)
    (1.0657,0.3159)
    (1.0667,0.3175)
    (1.0677,0.3191)
    (1.0687,0.3208)
    (1.0697,0.3226)
    (1.0707,0.3245)
    (1.0717,0.3266)
    (1.0727,0.3289)
    (1.0737,0.3313)
    (1.0747,0.3339)
    (1.0758,0.3367)
    (1.0768,0.3397)
    (1.0778,0.3428)
    (1.0788,0.3458)
    (1.0798,0.3487)
    (1.0808,0.3512)
    (1.0818,0.3535)
    (1.0828,0.3555)
    (1.0838,0.3573)
    (1.0848,0.3588)
    (1.0859,0.3603)
    (1.0869,0.3617)
    (1.0879,0.3630)
    (1.0889,0.3643)
    (1.0899,0.3655)
    (1.0909,0.3667)
    (1.0919,0.3677)
    (1.0929,0.3687)
    (1.0939,0.3697)
    (1.0949,0.3707)
    (1.0960,0.3718)
    (1.0970,0.3730)
    (1.0980,0.3743)
    (1.0990,0.3756)
    (1.1000,0.3770)
};
\addlegendentry{\texttt{CGN}}
\end{axis}
\end{tikzpicture}
\caption{Performance profiles for the four solver families with respect to solution quality.}
\label{fig:performance-profile}
\end{figure}

In this section, we compare solver families by taking for each instance the best solution obtained within each family. Figure~\ref{fig:performance-profile} presents the performance profile of the four solver families on all 972 instances. \texttt{PTN-DC} is the best-performing family overall, followed by \texttt{DPTN}, then \texttt{PTN}, and finally \texttt{CGN}. \texttt{PTN-DC} finds the best solution in over 83\% of the instances, and always finds a solution within 5\% of the best known solution. \texttt{CGN}'s performance is inconsistent: it finds the best solution in 8\% of instances, but also fails to find any solution at all in over 50\% of instances, with column generation often not even able to solve the root node of the model.

Figure~\ref{fig:network-performance-profiles} presents the performance profile per network. \texttt{PTN-DC} is the dominant model across all networks. \texttt{DPTN}'s performance is relatively stronger on Netherlands-IC-SPR and Switzerland-Large, the two largest instances, demonstrating that it also scales well. \texttt{CGN} completely fails on the largest instances, but is competitive on Randstad-IC, the smallest network, where it even finds 37\% of the best solutions. These results can be explained using the sizes of the PTN, DPTN and CGN, given in Table~\ref{tab:network-sizes}. Even with the smallest line pools, the CGN is much larger than the PTN or DPTN, such that the resulting model is much harder to solve. 
\input{figs/profilePerNetwork}
\begin{table}[t]
\centering
\caption{Network sizes. $|\mathcal{S}|$/$|\mathcal{C}|$: nodes/arcs in the PTN; $|\mathcal{W}|$/$|\mathcal{E}|$: nodes/arcs of the DPTN; $|\mathcal{V}|$/$|\mathcal{A}|$: nodes/arcs of the CGN for each line-pool size (S/M/L).}
\label{tab:network-sizes}
\resizebox{\textwidth}{!}{%
\begin{tabular}{lrr rr rr rr rr}
\toprule
& \multicolumn{2}{c}{PTN} & \multicolumn{2}{c}{DPTN} & \multicolumn{2}{c}{CGN-S} & \multicolumn{2}{c}{CGN-M} & \multicolumn{2}{c}{CGN-L} \\
\cmidrule(lr){2-3} \cmidrule(lr){4-5} \cmidrule(lr){6-7} \cmidrule(lr){8-9} \cmidrule(lr){10-11}
Network & $|\mathcal{S}|$ & $|\mathcal{C}|$ & $|\mathcal{W}|$ & $|\mathcal{E}|$ & $|\mathcal{V}|$ & $|\mathcal{A}|$ & $|\mathcal{V}|$ & $|\mathcal{A}|$ & $|\mathcal{V}|$ & $|\mathcal{A}|$ \\
\midrule
Randstad-IC          &   29 &   80 &    218 &     604 &    995 &    3,610 &   1,741 &    6,870 &   2,079 &    8,340 \\
Randstad-IC-SPR      &  109 &  324 &    866 &    2,384 &   5,293 &   21,870 &  12,905 &   57,610 &  18,753 &   85,510 \\
Netherlands-IC       &   78 &  196 &    548 &    1,422 &  10,646 &   46,160 &  29,828 &  135,310 &  66,878 &  309,580 \\
Netherlands-IC-SPR   &  256 &  750 &   2,012 &    5,366 &  18,128 &   77,100 &  55,408 &  251,550 & 125,846 &  585,580 \\
\midrule
Switzerland-Small        &   20 &   80 &    200 &     770 &   1,574 &    5,900 &   2,678 &   10,570 &   3,952 &   16,070 \\
Switzerland-Large     &   93 &  302 &    790 &    2,508 &  32,555 &  149,370 & 180,237 &  846,670 & 460,693 & 2,178,480 \\
\bottomrule
\end{tabular}}
\end{table}

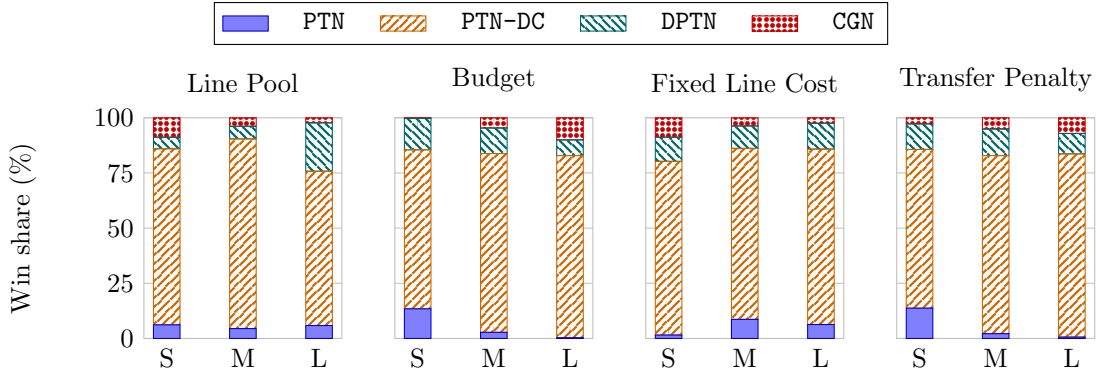
\begin{figure}[t]
\centering
\begin{tikzpicture}
\begin{axis}[hide axis, scale only axis, width=0.1pt, height=0.1pt,
    xmin=0, xmax=1, ymin=0, ymax=1,
    legend columns=4,
    legend style={column sep=1em, at={(0.5,0)}, anchor=north, font=\small},
    legend cell align={left}]
\addlegendimage{area legend, fill=blue!50!white, draw=blue!80!black}
\addlegendentry{\texttt{PTN}}
\addlegendimage{area legend, fill=orange!50!white, pattern={Lines[angle=45,  distance=2.3pt, line width=0.8pt]}, pattern color=orange!80!black, draw=orange!80!black}
\addlegendentry{\texttt{PTN-DC}}
\addlegendimage{area legend, fill=teal!50!white, pattern={Lines[angle=-45, distance=2.3pt, line width=0.8pt]}, pattern color=teal!80!black, draw=teal!80!black}
\addlegendentry{\texttt{DPTN}}
\addlegendimage{area legend, fill=red!50!white, pattern={Dots[distance=2.3pt, radius=0.8pt]}, pattern color=red!80!black, draw=red!80!black}
\addlegendentry{\texttt{CGN}}
\end{axis}
\end{tikzpicture}
\\[4pt]
\begin{tikzpicture}
\begin{groupplot}[
    group style={group size=4 by 1, horizontal sep=0.7cm},
    ybar stacked,
    width=4.2cm,
    height=4.5cm,
    ymin=0, ymax=100,
    ytick={0,25,50,75,100},
    yticklabel style={font=\small},
    xticklabel style={font=\small},
    axis line style={gray!50},
    tick style={gray!50},
    enlarge x limits=0.15,
    clip=false,
]
\nextgroupplot[title={\small Line Pool}, symbolic x coords={S,M,L}, xtick=data]
\addplot[fill=blue!50!white, draw=blue!80!black] coordinates {(S,6.1728) (M,4.4753) (L,5.8642)};
\addplot[fill=orange!50!white, pattern={Lines[angle=45,  distance=2.3pt, line width=0.8pt]}, pattern color=orange!80!black, draw=orange!80!black] coordinates {(S,79.7840) (M,85.9568) (L,69.9074)};
\addplot[fill=teal!50!white, pattern={Lines[angle=-45, distance=2.3pt, line width=0.8pt]}, pattern color=teal!80!black, draw=teal!80!black] coordinates {(S,5.2469) (M,5.7099) (L,21.9136)};
\addplot[fill=red!50!white, pattern={Dots[distance=2.3pt, radius=0.8pt]}, pattern color=red!80!black, draw=red!80!black] coordinates {(S,8.7963) (M,3.8580) (L,2.3148)};
\nextgroupplot[title={\small Budget}, symbolic x coords={S,M,L}, xtick=data, yticklabels={}]
\addplot[fill=blue!50!white, draw=blue!80!black] coordinates {(S,13.4259) (M,2.7778) (L,0.3086)};
\addplot[fill=orange!50!white, pattern={Lines[angle=45,  distance=2.3pt, line width=0.8pt]}, pattern color=orange!80!black, draw=orange!80!black] coordinates {(S,72.0679) (M,81.0185) (L,82.5617)};
\addplot[fill=teal!50!white, pattern={Lines[angle=-45, distance=2.3pt, line width=0.8pt]}, pattern color=teal!80!black, draw=teal!80!black] coordinates {(S,14.1975) (M,11.5741) (L,7.0988)};
\addplot[fill=red!50!white, pattern={Dots[distance=2.3pt, radius=0.8pt]}, pattern color=red!80!black, draw=red!80!black] coordinates {(S,0.3086) (M,4.6296) (L,10.0309)};
\nextgroupplot[title={\small Fixed Line Cost}, symbolic x coords={S,M,L}, xtick=data, yticklabels={}]
\addplot[fill=blue!50!white, draw=blue!80!black] coordinates {(S,1.5432) (M,8.6420) (L,6.3272)};
\addplot[fill=orange!50!white, pattern={Lines[angle=45,  distance=2.3pt, line width=0.8pt]}, pattern color=orange!80!black, draw=orange!80!black] coordinates {(S,78.7037) (M,77.4691) (L,79.4753)};
\addplot[fill=teal!50!white, pattern={Lines[angle=-45, distance=2.3pt, line width=0.8pt]}, pattern color=teal!80!black, draw=teal!80!black] coordinates {(S,10.9568) (M,10.1852) (L,11.7284)};
\addplot[fill=red!50!white, pattern={Dots[distance=2.3pt, radius=0.8pt]}, pattern color=red!80!black, draw=red!80!black] coordinates {(S,8.7963) (M,3.7037) (L,2.4691)};
\nextgroupplot[title={\small Transfer Penalty}, symbolic x coords={S,M,L}, xtick=data, yticklabels={}]
\addplot[fill=blue!50!white, draw=blue!80!black] coordinates {(S,13.7346) (M,2.1605) (L,0.6173)};
\addplot[fill=orange!50!white, pattern={Lines[angle=45,  distance=2.3pt, line width=0.8pt]}, pattern color=orange!80!black, draw=orange!80!black] coordinates {(S,71.9136) (M,80.7099) (L,83.0247)};
\addplot[fill=teal!50!white, pattern={Lines[angle=-45, distance=2.3pt, line width=0.8pt]}, pattern color=teal!80!black, draw=teal!80!black] coordinates {(S,11.5741) (M,12.0370) (L,9.2593)};
\addplot[fill=red!50!white, pattern={Dots[distance=2.3pt, radius=0.8pt]}, pattern color=red!80!black, draw=red!80!black] coordinates {(S,2.7778) (M,5.0926) (L,7.0988)};
\end{groupplot}
\node[rotate=90, anchor=center] at
    ($(group c1r1.west) + (-1.6cm,0)$)
    {\small Win share (\%)};
\end{tikzpicture}
\caption{Win share per solver family by parameter level. Each bar shows the average fraction of instances won by each family (ties split equally).}
\label{fig:win-share-params}
\end{figure}

\texttt{CGN}'s performance is not robust, but depends heavily on the parameter settings. Figure~\ref{fig:win-share-params} presents the win share of each solver family per parameter value for line pool, budget, fixed line cost, and transfer penalty. \texttt{CGN} performs better on relatively easy instances, with a small line pool, high budget, or small fixed line cost. \texttt{PTN}, which ignores transfers entirely, performs well when the transfer penalty is small, but its performance deteriorates quickly at medium or large transfer penalties. \texttt{DPTN}, by contrast, becomes increasingly competitive as the line pool grows: its win share is highest on large line pools. We also examined performance profiles per OD
matrix setting (dense versus sparse), but this did not appear to influence the relative performance of the solvers.

% Preamble: \usepgfplotslibrary{groupplots}  \usetikzlibrary{calc,patterns,patterns.meta}  \pgfplotsset{compat=newest}

%\clearpage
\subsection{Detailed Solver Analysis}
\label{sec:detailed}

In this section, we analyze the relative performance of the solvers within each family. Figure~\ref{fig:solver-performance-profiles} presents the performance profile of every solver variant, measured against the best solution found across all solvers. Table~\ref{tab:ptndc-gaps} reports the median optimality gap per solver and per network. Two caveats on the interpretation of these gaps are in order. First, for solvers that rely on a priori route enumeration, the gap is computed relative to the restricted model with that route set. Since the route set may exclude good solutions, the reported gap can be optimistic even when it appears small. For the branch-and-price solvers on the DPTN and CGN, routes are generated on the fly and the LP relaxation is solved to optimality at each node, so their lower bounds are valid for their respective models. However, because only the CGN models transfers exactly, \texttt{CGN} is the sole solver whose lower and upper bounds are both valid for the true problem. Second, a small gap need not signal a good line plan: a restricted or approximate model is solving a proxy for the true problem, so its objective estimate may diverge from the true cost in either direction. We quantify the magnitude and sign of these estimation errors per solver in Table~\ref{tab:objective-estimate-error}.
\input{figs/profilePerSolver}

Comparing the small, medium, and large route sets for all solvers that require a priori route enumeration, the medium-sized route sets consistently perform well. Small route sets exclude too many useful routes, while large route sets tend to perform marginally worse than medium ones, likely because they introduce many redundant routes that complicate the solution process. For \texttt{PTN-DC} and \texttt{DPTN}, large route sets exhibit a modest advantage for small performance ratios, indicating that larger route sets are beneficial whenever the resulting model remains tractable. However, on larger instances the route sets become prohibitively large, and medium-sized route sets outperform them. This is confirmed by the optimality gaps: particularly on Netherlands-IC-SPR, the instance with the largest number of OD pairs, large route sets substantially increase gaps. On some instances with a large route set and a large line pool, constructing the basic direct connection constraints for the \texttt{PTN-DC} solvers is itself very time-consuming, leaving little or no time for actual optimization. This underscores the importance of tailoring the route set to the instance at hand.

Within the \texttt{PTN-DC} family, \texttt{PTN-FDC} --- the Benders decomposition model --- has a slight edge over \texttt{PTN-BDC}, which includes only the basic direct connection constraints. The median optimality gaps of the Benders approach are also substantially lower. That said, \texttt{PTN-BDC} still performs very well and may suffice for many practical applications.

Within the DPTN family, branch-and-price outperforms route enumeration with the small route set, but is itself outperformed by route enumeration with medium or large route sets. This resolves the question raised in
Section~\ref{sec:bnp}: for the DPTN, a priori route enumeration is preferable to branch-and-price, provided the route set is chosen well. The addition of strong inequalities clearly improves performance throughout. For \texttt{CGN}, strong inequalities substantially strengthen the model, raising the fraction of instances for which \texttt{CGN} finds the best known solution to 8\% and reducing optimality gaps. On harder instances, however, they also make it significantly more difficult to find any feasible solution, as evidenced by the sharp increase in infinite gaps in Table~\ref{tab:ptndc-gaps} for \texttt{CGN-BP\&C} relative to \texttt{CGN-B\&P}. Intuitively, strong inequalities are always desirable in principle, but their number grows with the number of arcs in the routing graph. Since the CGN contains far more arcs than the DPTN, separating and managing cuts becomes a bottleneck on large instances, offsetting the benefit of the tighter relaxation. For such instances, a different solver family would be preferable in any case.

% Preamble: \usepackage{booktabs}

\begin{table}[t]
\centering
\resizebox{\textwidth}{!}{%
\begin{tabular}{lrrrrrrrrrrrrrrrr}
\toprule
 & \multicolumn{3}{c}{\texttt{PTN}} & \multicolumn{3}{c}{\texttt{PTN-BDC}} & \multicolumn{3}{c}{\texttt{PTN-FDC}} & \multicolumn{5}{c}{\texttt{DPTN}} & \multicolumn{2}{c}{\texttt{CGN}} \\
\cmidrule(lr){2-4} \cmidrule(lr){5-7} \cmidrule(lr){8-10} \cmidrule(lr){11-15} \cmidrule(lr){16-17}
Network & \texttt{S} & \texttt{M} & \texttt{L} & \texttt{S} & \texttt{M} & \texttt{L} & \texttt{S} & \texttt{M} & \texttt{L} & \texttt{S} & \texttt{M} & \texttt{L} & \texttt{B\&P} & \texttt{BP\&C} & \texttt{B\&P} & \texttt{BP\&C} \\
\midrule
Randstad-IC & 0.00 & 0.00 & 0.00 & 0.00 & 0.00 & 0.00 & 0.00 & 0.00 & 0.00 & 0.00 & 0.00 & 0.00 & 2.50 & 0.27 & 12.75 & 1.73 \\Randstad-IC-SPR & 0.00 & 0.00 & 0.00 & 0.00 & 0.11 & 0.27 & 0.00 & 0.10 & 0.26 & 0.00 & 0.01 & 0.05 & 2.17 & 1.16 & 12.86 & $\infty$ \\Netherlands-IC & 0.00 & 0.18 & 0.15 & 0.04 & 1.92 & 2.05 & 0.18 & 1.21 & 1.05 & 0.01 & 1.35 & 1.11 & 6.83 & 4.47 & $\infty$ & $\infty$ \\Netherlands-IC-SPR & 0.03 & 0.10 & 1.76 & 0.07 & 2.03 & 52.92 & 0.10 & 1.00 & 21.21 & 0.04 & 1.05 & 14.36 & 3.01 & $\infty$ & $\infty$ & $\infty$ \\Switzerland-Small & 0.00 & 0.16 & 0.07 & 0.00 & 0.45 & 0.39 & 0.00 & 0.07 & 0.06 & 0.00 & 0.19 & 0.20 & 1.52 & 0.27 & 8.89 & 5.14 \\Switzerland-Large & 0.00 & 0.25 & 0.14 & 0.17 & 2.41 & 1.46 & 0.19 & 0.97 & 0.61 & 0.00 & 1.12 & 0.77 & 12.42 & 0.96 & $\infty$ & $\infty$ \\
\bottomrule
\end{tabular}%
}
\caption{Median optimality gap (\%) for \texttt{PTN}, \texttt{PTN$-$DC}, and \texttt{DPTN} solver variants per network. Infeasible runs and runs without a valid lower bound are treated as $\infty$.}
\label{tab:ptndc-gaps}
\end{table}

To provide insight into which problem characteristics complicate finding 
high-quality solutions, Figure~\ref{fig:gap-param-profiles} shows how the 
gap profile depends on the line pool, budget, fixed line cost, and transfer penalty settings, using \texttt{PTN-FDC-M} as the solver. Larger line pools and smaller budgets make instances progressively harder to solve, as they increase the number of candidate lines and tighten the feasible region, respectively. The effect of fixed line cost is more abrupt: instances with no fixed cost are considerably easier, while medium and large fixed costs yield similar gaps. This is consistent with the fact that fixed costs introduce binary activation variables that complicate the MIP, regardless of their magnitude. A similar threshold pattern holds for the transfer penalty: a small penalty is easier to solve, likely because transfers are nearly free and passenger routing is less constrained, but there is little difference between medium and large penalties.

\input{figs/gapParameter}

Beyond solution quality and optimality gaps, it is instructive to examine how accurately each solver \textit{estimates} the objective value of the line plans it produces. Table~\ref{tab:objective-estimate-error} reports the median signed relative error $(\hat{z} - z^*) / z^* \times 100$ for each solver variant and network, where positive values indicate an overestimate and negative values an underestimate.

All route-enumeration solvers tend to overestimate with the small route set, since the restricted route set excludes good indirect routings and the \texttt{CGN} evaluation then finds cheaper solutions than the solver anticipated. This effect diminishes with larger route sets and, in some cases, reverses into an underestimate. \texttt{PTN} ignores transfers entirely, so the \texttt{CGN} evaluation reveals transfer penalties not accounted for during optimization. \texttt{PTN-DC} and \texttt{DPTN} partially capture transfers, so with medium or large route
sets their estimates are closer to the true objective. In particular, both variants of \texttt{PTN-DC} provide highly accurate estimates, with median errors typically below 1\% in absolute value.
% Preamble: \usepackage{booktabs}

\begin{table}[t]
\centering
\resizebox{\textwidth}{!}{%
\begin{tabular}{lrrrrrrrrrrrrrrrr}
\toprule
 & \multicolumn{3}{c}{\texttt{PTN}} & \multicolumn{3}{c}{\texttt{PTN-BDC}} & \multicolumn{3}{c}{\texttt{PTN-FDC}} & \multicolumn{5}{c}{\texttt{DPTN}} & \multicolumn{2}{c}{\texttt{CGN}} \\
\cmidrule(lr){2-4} \cmidrule(lr){5-7} \cmidrule(lr){8-10} \cmidrule(lr){11-15} \cmidrule(lr){16-17}
Network & \texttt{S} & \texttt{M} & \texttt{L} & \texttt{S} & \texttt{M} & \texttt{L} & \texttt{S} & \texttt{M} & \texttt{L} & \texttt{S} & \texttt{M} & \texttt{L} & B\&P & BP\&C & B\&P & BP\&C \\
\midrule
Randstad-IC & +8.5 & -9.9 & -9.9 & +19.2 & +0.1 & -0.0 & +19.2 & +0.1 & -0.1 & +18.2 & -2.7 & -4.0 & -3.2 & -4.1 & +0.0 & +0.0 \\Randstad-IC-SPR & +21.8 & -7.8 & -8.5 & +38.2 & +0.6 & -0.4 & +38.3 & +0.6 & -0.6 & +34.0 & -2.6 & -4.3 & -3.8 & -4.5 & +0.0 & -- \\Netherlands-IC & +26.0 & -7.0 & -7.3 & +34.5 & -0.0 & -0.2 & +34.8 & -0.3 & -0.4 & +34.9 & -3.8 & -4.0 & -2.2 & -3.8 & +0.0 & -- \\Netherlands-IC-SPR & +16.1 & -4.4 & -7.3 & +30.9 & +0.8 & -0.1 & +31.1 & +1.0 & -0.8 & +25.8 & -1.1 & -2.9 & -4.2 & -4.2 & -- & -- \\Switzerland-Small & +19.8 & -2.7 & -3.6 & +22.6 & +0.1 & -0.2 & +22.2 & -0.2 & -0.4 & +24.3 & -0.3 & -1.6 & -1.7 & -2.0 & +0.0 & +0.0 \\Switzerland-Large & +29.1 & -5.0 & -7.1 & +36.9 & +1.5 & +0.1 & +36.8 & +1.4 & +0.0 & +41.1 & +1.4 & -3.1 & +0.1 & -4.6 & -- & -- \\
\bottomrule
\end{tabular}%
}
\caption{Median signed relative error (\%) of the objective estimate: $(\hat{z} - z^*) / z^* \times 100$. Negative values indicate an underestimate (optimistic); positive values an overestimate. Runs without a feasible solution or missing estimate are excluded.}
\label{tab:objective-estimate-error}
\end{table}

\subsection{Impact of Demand Sparsification}
\label{sec:demSpar}
This section analyzes the impact of the preprocessing techniques of 
Section~\ref{sec:preprocessing}. We use \texttt{PTN-FDC-M} as the solver. Figure~\ref{fig:cross-eval-gap-profiles} shows optimality gaps for the 
original and sparsified demand matrices. The sparsified matrix yields 
consistently smaller or comparable gaps across most networks, with the 
most pronounced improvement on Netherlands-IC-SPR, where the reduced 
number of OD pairs makes the instance substantially easier to solve.
\input{figs/gapODs}

To assess whether better gaps on the sparsified instance translate into 
better solutions to the original problem, we re-evaluate all line plans 
on the full OD matrix. Figure~\ref{fig:cross-eval-network-profiles} shows 
the resulting performance profiles. For most networks, solving on the 
original matrix performs better, since these 
instances are already tractable. For Netherlands-IC-SPR, however, the 
line plans obtained on the sparsified matrix are substantially better, 
demonstrating that preprocessing the OD matrix can yield considerably 
improved solutions on large instances where the full matrix is a 
computational bottleneck. 
\input{figs/profileODs}
\section{Conclusion}
\label{sec:con}

This paper provided a theoretical and empirical comparison of the main 
network-based models for line planning with passenger routing. On the 
theoretical side, we made the structural relationships between the models 
explicit and established a hierarchy in terms of solution quality. On the 
empirical side, we evaluated all models on 972 instances drawn from the Dutch and Swiss railway networks.

Our results offer concrete guidance for line planning researchers and practitioners. The PTN-based direct connection model is the most robust and best-performing approach across the full range of instances. The exact model based on the change-and-go network is competitive only on small instances and frequently fails to find any feasible solution on larger ones. These findings suggest that the additional complexity of the CGN is rarely warranted in practice, and that carefully designed PTN-based approximations are preferable for line planning at scale. We therefore believe the PTN-based approximations are the most promising foundation for future extensions of the line planning problem.
\bibliographystyle{abbrvnat}
\bibliography{references}

\end{document}